\documentclass[11pt,a4paper,leqno]{amsart}
\usepackage{amssymb}
\numberwithin{equation}{section}

\newtheoremstyle{theorem}{3pt}{3pt}%
{\it}%         Body font
{}%         Indent amount (empty = no indent, \parindent = para indent)
{\bfseries}% Thm head font  (could be also \sc)
{:}%        Punctuation after thm head
{.5em}%     Space after thm head (\newline = linebreak)
{}%         Thm head spec
\theoremstyle{theorem}
\newtheorem{theorem}{Theorem}[section]

\newtheorem{proposition}[theorem]{Proposition}

\newtheorem{definition}[theorem]{Definition}

\newtheoremstyle{example}{3pt}{3pt}%
{}%         Body font
{}%         Indent amount (empty = no indent, \parindent = para indent)
{\sc}% Thm head font  (could be also \sc)
{:}%        Punctuation after thm head
{.5em}%     Space after thm head (\newline = linebreak)
{}%         Thm head spec
\theoremstyle{example}
\newtheorem{example}[theorem]{Example}
\newtheoremstyle{remark}{3pt}{3pt}%
{}%         Body font
{}%         Indent amount (empty = no indent, \parindent = para indent)
{\sc}% Thm head font  (could be also \sc)
{:}%        Punctuation after thm head
{.5em}%     Space after thm head (\newline = linebreak)
{}%         Thm head spec
\theoremstyle{remark}
\newtheorem{remark}{Remark}[section]

\numberwithin{equation}{section}

\newcommand{\acknowledge}{\subsection*{Acknowledgments}}
\newcommand{\thismonth}{\ifcase\month\or
  January\or February\or March\or April\or May\or June\or
  July\or August\or September\or October\or November\or December\fi
  \space\number\year}
\makeatletter
\newcommand{\low}{\@ifnextchar^{}{^{\vphantom x}}}
\newcommand{\high}{\@ifnextchar_{}{_{\vphantom I}}}
\makeatother
\DeclareSymbolFont{script}{U}{eus}{m}{n}
\DeclareSymbolFontAlphabet{\mathscr}{script}
\DeclareMathSymbol{\EuWedge}{0}{script}{"5E}
\DeclareMathAlphabet{\mathrmsl}{OT1}{cmr}{m}{sl}
\newcommand{\rssymb}[2]{\newcommand{#1}{{\mathrmsl{#2}}}}
\newcommand{\calsymb}[2]{\newcommand{#1}{{\mathcal{#2}}}}
\newcommand{\bbsymb}[2]{\newcommand{#1}{{\mathbb{#2}}}}
\newcommand{\lieoper}[2]{\newcommand{#1}{\mathop
  {\mathfrak{#2}\null}\nolimits}}
\newcommand{\oper}[3][n]{\newcommand{#2}{\mathop
  {\mathrm{#3}\null}\ifx n#1\nolimits\else\limits\fi}}
\newcommand{\rsoper}[3][n]{\newcommand{#2}{\mathop
  {\mathrmsl{#3}\null}\ifx n#1\nolimits\else\limits\fi}}
\bbsymb\C{C} \bbsymb\F{F} \bbsymb\HQ{H}\bbsymb\N{N} \bbsymb\Q{Q}
\bbsymb\R{R} \bbsymb\U{U} \bbsymb\V{V} \bbsymb\W{W} \bbsymb\Z{Z}
\bbsymb\bbf{F} \bbsymb\bbk{K} \bbsymb\bbi{I} \bbsymb\bbl{L} \bbsymb\bbo{O}
\bbsymb\bbj{J} 
\bbsymb\bby{Y} 
\bbsymb\bbp{P} 
\bbsymb\bba{A} 
\calsymb\cA{A} \calsymb\cB{B} \calsymb\cC{C} \calsymb\cD{D} \calsymb\cE{E}
\calsymb\cF{F} \calsymb\cG{G} \calsymb\cH{H} \calsymb\cI{I} \calsymb\cJ{J}
\calsymb\cK{K} \calsymb\cL{L} \calsymb\cM{M} \calsymb\cN{N} \calsymb\cO{O}
\calsymb\cP{P} \calsymb\cQ{Q} \calsymb\cR{R} \calsymb\cS{S} \calsymb\cT{T}
\calsymb\cU{U} \calsymb\cV{V} \calsymb\cW{W} \calsymb\cX{X} \calsymb\cY{Y}
\calsymb\cZ{Z}
\newcommand{\eps}{\varepsilon}

\renewcommand{\geq}{\geqslant} \renewcommand{\leq}{\leqslant}
\oper\End{End} \oper\Hom{Hom}                    % Vector space constructions
\oper\Sym{Sym} \oper\Skew{Skew}
\oper\Aut{Aut}                                   % Group constructions
\oper\GL{GL} \oper\SL{SL}\oper\Symp{Sp}
\oper\CO{CO} \oper\On{O} \oper\SO{SO} \oper\Pin{Pin} \oper\Spin{Spin}
\oper\CU{CU} \oper\Un{U} \oper\SU{SU} \oper\PSU{PSU}
\rsoper\Diff{Diff} \rsoper\SDiff{SDiff}
\lieoper\der{der}                                % Lie algebra constructions
\lieoper\gl{gl} \lieoper\sgl{sl}\lieoper\symp{sp}
\lieoper\co{co} \lieoper\so{so} \lieoper\spin{spin}
\lieoper\cu{cu} \lieoper\un{u}  \lieoper\su{su}
\rsoper\Vect{Vect} \rsoper\Ham{Ham}
\def\la#1{\hbox to #1pc{\leftarrowfill}}
\def\ra#1{\hbox to #1pc{\rightarrowfill}}

\newcommand{\ip}[1]{\langle#1\rangle}

\newcommand{\norm}[2][]{|\mkern-2mu|#2|\mkern-2mu|
  _{\lower1pt\hbox{${}_{#1}$}}}
\newcommand{\Norm}[2][]{\bigl|\mkern-3mu\bigr|#2\bigr|\mkern-3mu\bigr|
  _{\lower1pt\hbox{${}_{#1}$}}}

\newcommand{\dual}{^{*\!}}

\newcommand{\dsum}{\oplus}                  % small direct sum
\newcommand{\tens}{\otimes}                 % small tensor product
\newcommand{\Proj}{\mathrmsl{P}}            % projective
\newcommand{\CP}[1]{\C\Proj^{#1}}           % complex projective space
\newcommand{\HP}[1]{\HQ\Proj^{#1}}          % quaternionic projective space
\newcommand{\RH}[1]{\cH^{#1}}               % real hyperbolic space
\newcommand{\CH}[1]{\C\cH^{#1}}             % complex hyperbolic space
\newcommand{\HH}[1]{\HQ\cH^{#1}}            % quaternionic hyperbolic space
\newcommand{\Cinf}{\mathrm{C}^\infty}       % smooth sections
\rsoper\dimn{dim}                           % dimension
\rsoper\grad{grad}                          % gradient
\rsoper\kernel{ker}\rsoper\image{im}        % kernel and image
\rsoper\alt{alt}   \rsoper\sym{sym}         % alternating and symmetric part
\rsoper\Ad{Ad}     \rsoper\ad{ad}           % adjoint action or bundle
\rsoper\CoAd{CoAd} \rsoper\coad{coad}       % coadjoint action
\rsoper\trace{tr}  \rsoper\trfree{tf}       % trace and tracefree part
\rsoper\detm{det}                           % determinant
\rsoper\Vol{Vol}                            % volume
\rsoper\divg{div}                           % divergence
\rsoper\sign{sign}                          % sign function
\rssymb\iden{id}                            % identity
\rssymb\vol{vol}                            % volume element
\oper\Imag{Im}\oper\Real{Re}                % real and imaginary
\newcommand{\sd}{{\raise1pt\hbox{$\scriptscriptstyle +$}}}
\newcommand{\asd}{{\raise1pt\hbox{$\scriptscriptstyle -$}}}
\newcommand{\sdasd}{{\raise1pt\hbox{$\scriptscriptstyle\pm$}}}
\newcommand{\asdsd}{{\raise1pt\hbox{$\scriptscriptstyle\mp$}}}

\rsoper\scal{scal}
\def\kahl/{k\"ahler}
\def\Kahl/{K{\"a}hler}
\newcommand{\g}{{\mathfrak g}}
\newcommand{\bfp}{{\mathbf p}}
\newcommand{\bfu}{{\mathbf u}}
\newcommand{\bfv}{{\mathbf v}}
\newcommand{\bfw}{{\mathbf w}}
\newcommand{\bfx}{{\mathbf x}}
\newcommand{\bfy}{{\mathbf y}}
\newcommand{\bfz}{{\mathbf z}}
\newcommand{\bfa}{{\mathbf a}}
\newcommand{\bfb}{{\mathbf b}}
\def\decdnar#1{\phantom{\hbox{$\scriptstyle{#1}$}}
\left\downarrow\vbox{\vskip15pt\hbox{$\scriptstyle{#1}$}}\right.}
\begin{document}
\title[Toric self-dual Einstein metrics]
{Toric self-dual Einstein metrics as quotients}
\author[C. Boyer]{Charles P. Boyer}
\address{CPB and KG: Department of Mathematics and Statistics,
University of New Mexico,
Albuquerque, NM 87131, USA.}
\email{cboyer@math.unm.edu}
\email{galicki@math.unm.edu}
\author[D. Calderbank]{David M. J. Calderbank}
\address{DMJC: School of Mathematics\\
University of Edinburgh\\ King's Buildings, Mayfield Road\\
Edinburgh EH9 3JZ\\ Scotland.}
\email{davidmjc@maths.ed.ac.uk}
\author[K. Galicki]{Krzysztof Galicki}
\author[P. Piccinni]{Paolo Piccinni}
\address{PP: Universit\`a degli Studi di Roma, ``La Sapienza'', Dipartimento di
Matematica, Piazzale Aldo Moro, I-00185 Roma, Italia.}
\email{piccinni@mat.uniroma1.it}
\date{\thismonth}
\thanks{During the preparation of this work the first and third authors were
supported by NSF grant DMS-0203219. The second author was supported by the
Leverhulme Trust, the William Gordon Seggie Brown Trust and an EPSRC Advanced
Fellowship. The fourth author was supported by the MIUR Project ``Propriet\`a
Geometriche delle Variet\`a Reali e Complesse''.  The authors are also grateful
for support from EDGE, Research Training Network HPRN-CT-2000-00101, funded by
the European Human Potential Programme.}
\begin{abstract}
We use the quaternion K\"ahler reduction technique to study old and new
self-dual Einstein metrics of negative scalar curvature with at least a
two-dimensional isometry group, and relate the quotient construction to the
hyperbolic eigenfunction Ansatz. We focus in particular on the
(semi-)quaternion K\"ahler quotients of (semi-)quaternion K\"ahler
hyperboloids, analysing the completeness and topology, and relating them to
the self-dual Einstein Hermitian metrics of Apostolov--Gauduchon and Bryant.
\end{abstract}
\maketitle
\vspace{-2mm}
\section*{Introduction}

There has been quite a lot of interest in self-dual Einstein (SDE) metrics in
dimension four. In the negative scalar curvature case, such metrics naturally
generalize the symmetric metrics on the real $4$-ball $\RH4\simeq\HH1$ (the
real or quaternionic hyperbolic metric) and on the complex $2$-disc $\CH2$
(the complex hyperbolic or Bergman metric).

A rather general construction of negative SDE metrics was offered by C. LeBrun
in 1982 \cite{MR83d:83019}. LeBrun observed that for any real-analytic
conformal structure $[h]$ on $S^3$, there is a Riemannian metric $g_0$ defined
on some open neighborhood of $S^3\subset\R^4$ such that $g_0$ is self-dual,
the restriction of it to $S^3$ is in the conformal class $[h]$, and moreover
$g=f^{-2}g_0$ is Einstein for some defining function $f$ for $S^3$ in this open
neighbourhood. However, this result is purely local: the Einstein metric it
defines typically cannot be extended to a complete metric everywhere inside
the ball.

Nevertheless, in later work \cite{MR92i:53042}, LeBrun showed that the moduli
space of negative complete SDE metrics on a ball is infinite dimensional,
which led him to formulate a conjecture. A conformal structure on $S^3$ is
said to have positive frequency if it bounds a complete SDE metric on the ball
and negative frequency if bounds a complete anti-self-dual Einstein (ASDE)
metric on the ball. The conjecture then asserts that near the standard
conformal structure on $S^3$ (in an appropriate sense) the moduli spaces of
positive and negative frequency subspaces are transverse (i.e., their tangent
spaces at the standard conformal structure give a direct sum
decomposition). The positive frequency conjecture is now proven, thanks to the
remarkable work of O. Biquard \cite{MR1908060} (see also \cite{MR2001k:53079,
MR2001k:53080}).  However, this still provides very little information about
which conformal structures on $S^3$ bound complete SDE metrics on the ball.

The known examples are rather few. Apart from the hyperbolic metric, the first
such metrics were obtained by H. Pedersen in \cite{MR87i:53070}: the conformal
class $[h]$ on $S^3$ is represented by a Berger sphere metric
$\sigma_1^2+\sigma_2^2+\lambda^2\sigma_3^2$ (where $\sigma_1$, $\sigma_2$ and
$\sigma_3$ are the standard left-invariant $1$-forms on $S^3\simeq \Symp(1)$
and $\lambda$ is a nonzero constant), and the corresponding complete SDE
metric on the $4$-ball is equally explicit.  Later, N. Hitchin
\cite{MR96g:53057} generalized this result by showing that any left-invariant
conformal structure on $S^3$ determines a complete SDE metric on the ball,
although now explicitness requires elliptic rather than elementary functions.

The reason that these metrics are tractible is the presence of symmetry.  The
real and complex hyperbolic metrics have isometry groups $SO(4,1)$ and
$U(2,1)$ respectively, while the Pedersen metrics on the ball have isometry
group $U(2)$. There are also (related) $U(2)$-invariant SDE metrics on complex
line bundles $\cO(n)\rightarrow\CP2$, with $n\geq3$, called the
Pedersen--LeBrun metrics~\cite{MR89f:53107, MR96g:53057}. Hitchin
\cite{MR96g:53057} actually classifies all $\SU(2)$-invariant SDE metrics, and
proves that the complete examples of negative scalar curvature consist only of
the real and complex hyperbolic metrics, the Pedersen and Pedersen--LeBrun
metrics, and SDE metrics on the ball associated to a left-invariant conformal
or CR structure on $S^3\simeq SU(2)$.

Recent progress on SDE metrics with symmetry concerns the much smaller
symmetry group $T^2\simeq S^1\times S^1$ (and its non-compact forms): such SDE
metrics are said to be \emph{toric}. In the positive case, these metrics can
be constructed using the Galicki--Lawson quaternion K\"ahler reduction
\cite{MR89m:53075} of the quaternionic projective space $\HQ\Proj^n$ by the
action of an $(n-1)$-dimensional subtorus of the maximal torus of
$\Symp(n+1)$.  Although the only positive SDE metrics on compact manifolds are
the standard metrics on $S^4$ and $\CP2$, these methods produce positive SDE
metrics on compact orbifolds.  The general such metrics were described by
C. Boyer et al. in \cite{MR99b:53066}, following the construction by Galicki
and Lawson of positive SDE metrics on weighted projective spaces
\cite{MR89m:53075}.  It is natural to conjecture that all positive compact
SDE orbifolds arise in this way: this would be similar to a related result of
R. Bielawski \cite{MR2001c:53058} stating that all toric 3-Sasakian manifolds
(in any dimension) are the 3-Sasakian quotients considered in
\cite{MR99b:53066}.

Another impetus to study toric SDE metrics comes from the recent work
\cite{MR1950174} of D. Calderbank and H. Pedersen, who proved that if a
(positive or negative) SDE metric admits two commuting Killing vector fields,
it can be expressed locally in an explicit form depending on a single function
$F$ on the upper-half plane, where $F$ is an eigenfunction of the hyperbolic
Laplacian with eigenvalue $3/4$. Conversely, any metric of this form is an SDE
metric. Calderbank and Pedersen then showed explicitly how the positive SDE
metrics of Galicki--Lawson and Boyer et al. arise from such an eigenfunction
$F$, and tied together a number of examples of negative SDE metrics.

The (locally) toric SDE metrics of \cite{MR1950174} also relate to a recent
study by V. Apostolov and P. Gauduchon of SDE Hermitian metrics
\cite{math.DG/0003162}. SDE metrics with symmetry are conformal to metrics
which are K\"ahler with the opposite orientation (hence scalar-flat), but it
is much rarer for an SDE metric to admit a Hermitian structure inducing the
given orientation.  Nevertheless, many of the examples of SDE metrics
discussed so far are Hermitian in this sense. Other non-locally symmetric
examples of SDE Hermitian metrics include cohomogeneity one metrics under the
action of $\R\times {\rm Isom}(\R^2)$, $U(1,1)$, and $U(2)$ constructed by
A. Derdzi\'nski \cite{Der:em} (the $U(2)$ case being the Pedersen--LeBrun
metrics mentioned above). Apostolov and Gauduchon show, quite generally, that
SDE Hermitian metrics always admit two distinguished commuting Killing vector
fields, and that if the induced local $\R^2$ action does not have two
dimensional generic orbits, then the isometry group necessarily acts
transitively or with cohomogeneity one. In either case, they show that SDE
Hermitian metrics are toric, hence given locally by the metrics of Calderbank
and Pedersen.

The emergence of non-trivial isometries for SDE Hermitian metrics is perhaps
less surprising in view of a link with recent work of R. Bryant on
Bochner-flat K\"ahler metrics \cite{MR2002i:53096}. In four dimensions, the
Bochner tensor coincides with the anti-self-dual Weyl tensor and so K\"ahler
metrics with vanishing Bochner tensor are just self-dual K\"ahler metrics.
Apostolov and Gauduchon show that SDE Hermitian metrics are necessarily
conformal to self-dual K\"ahler metrics, hence they belong to the class of
metrics studied by Bryant. In his impressive paper, Bryant obtains an explicit
local classification of Bochner-flat K\"ahler metrics and studies in detail
their global geometry. The symmetries here arise naturally from a differential
system, which amounts to the realisation of Bochner-flat K\"ahler
$2n$-manifolds as local quotients of the flat CR structure on $S^{2n+1}$.
Bryant's work not only provides an alternative way of classifying SDE
Hermitian metrics locally, but it also gives insight into the question of
completeness, and he discusses some examples in an appendix to his paper.

In spite of this work (and in contrast to the case of $\SU(2)$ symmetry, where
Hitchin provides a classification) the issue of completeness for negative SDE
Hermitian metrics is not yet fully explored, and for the toric SDE metrics in
general, the complete examples are far from understood. In fact, there are
very many examples. In \cite{CaSi:emcs}, Calderbank and M. A. Singer
constructed examples of complete SDE metrics on resolutions of complex cyclic
singularities and showed that the moduli of such metrics is (continuously)
infinite dimensional. In particular these metrics can have arbitrarily large
second Betti number (cf. \cite{MR99b:53066} in the positive case).  Examples
of infinite topological type are also known.

The simplest examples in \cite{CaSi:emcs} are quaternion K\"ahler quotients of
$\HH m$ generalizing the Pedersen--LeBrun metrics on $\cO(n)$ ($n\geq 3$), and
may be viewed as negative analogues of the compact orbifold SDE metrics of
Galicki--Lawson and Boyer et al.

In fact many of the metrics discussed in this introduction occur as quaternion
K\"ahler quotients \cite{MR88f:53088, MR89m:53075}. For positive toric SDE
metrics, compact orbifold examples are well understood (as we have discussed).
For negative toric SDE metrics, many examples have been introduced as
quotients by Galicki \cite{MR88j:53076, MR92i:53040}, but the quotient
approach has not been thoroughly explored.  Our purpose in this work is to
develop systematically the quotient approach to toric SDE metrics, which has a
number of advantages. In addition to producing an abundance of examples
locally, the quotient approach provides more direct insight into the global
behaviour of such metrics (completeness or topology), as well as a systematic
way to organise these examples into families.

In this paper we set the initial stage for such a systematic study by
considering the toric SDE metrics arising as (semi-)quaternion K\"ahler
quotients of $8$-dimensional quaternionic hyperbolic space $\HH2$ and its
indefinite signature analogue $\HH{1,1}$ by a one dimensional group action. A
given reduction may be encoded by the adjoint orbits in $\symp(1,2)$ of the
generator of the action, which in turn may be classified using work of
Burgoyne and Cushman \cite{MR55:5761}. There are essentially four distinct
possible types of generator:
\begin{enumerate}
\item elements belonging to the Lie algebra a maximal torus;
\item elements in a Cartan subalgebra with exponential image $S^1\times\R$;
\item non-semisimple elements with two step nilpotent part;
\item non-semisimple elements with three step nilpotent part.
\end{enumerate}
The quaternion K\"ahler quotients by generators in the first two classes
correspond to the $3$-pole solutions discussed in \cite{MR1950174}, but we
present a detailed and self-contained analysis of the completeness and
topology of the quotient. The other two classes may be regarded as limiting
cases, but the geometry of the quotient is less well studied.

According to Apostolov and Gauduchon \cite{math.DG/0003162}, quaternion
K\"ahler quotients of $\HH2$ (and $\HP2$) by one dimensional group actions are
SDE Hermitian, and their argument applies also to quotients of $\HH{1,1}$.
Therefore all of the quotients we discuss in this paper are SDE Hermitian
manifolds. Furthermore, by comparing our examples with the classification of
self-dual K\"ahler metrics by Bryant \cite{MR2002i:53096}, we see that in fact
all SDE Hermitian metrics with nonzero scalar curvature are (at least locally)
quaternion K\"ahler quotients of $\HP2$, $\HH2$ or $\HH{1,1}$.

In addition to studying the quotients of $\HH2$ and $\HH{1,1}$ in detail, we
develop some aspects of the general theory of quotients of $\HH{k,l}$ by
$(k+l-1)$-dimensional Abelian semi-quaternion K\"ahler group actions. In
particular we show how the quotient metrics are related to the hyperbolic
eigenfunction Ansatz, simplifying and extending a result of~\cite{MR1950174}.

\bigskip

\acknowledge The first author thanks the Ecole Polytechnique, Palaiseau and
the Universit\`a di Roma ``La Sapienza'' for hospitality and support. The
third author would like to thank the Universit\`a di Roma ``La Sapienza'',
I.N.d.A.M, M.P.I-Bonn, and IHES as parts of this paper were written during his
visits there. The fourth named author would like to thank University of New
Mexico for hospitality and support.  The authors are grateful to Paul
Gauduchon, Michael Singer and Pavel Winternitz for invaluable discussions.

\section{Semi-Quaternionic Projective Spaces}

\begin{definition}
Let $(M^{4n},g)$ be a semi-Riemannian manifold of signature $(4\nu,4n-4\nu)$.
We say that $(M^{4n},g)$ is semi-quaternion K\"ahler if the holonomy group of
the metric connection is a subgroup of $\Symp(\nu,n-\nu)\cdot \Symp(1)$ when
$n>1$. As usual, when $n=1$ we extend our definition and require that $(M,g)$
be self-dual and Einstein. We will always suppose that the scalar curvature of
$(M,g)$ is nonzero. We refer to $\nu$ as the quaternionic index of $M$.
\end{definition}

Exactly as in the Riemannian case, the above definition implies the existence
of the quaternion K\"ahler 4-form $\Omega$ which is parallel with respect to
the Levi-Civita connection and gives rise to the quaternionic rank 3 bundle
$\cV$ over $M$.

The simplest example of semi-quaternion K\"ahler manifolds are obtained as
follows. Let $\HQ^{k,l}=\{\bfu=(\bfa,\bfb)\ \ | \ \ \bfa
=(u_0,\ldots,u_{k-1}),\ \ \bfb=(u_k,\ldots,u_{k+l})\}$ be the set of all
quaternionic $(n+1)$-vectors together with the symmetric form
\begin{equation}\label{symform}
F_{k,l}(\bfu^1,\bfu^2)=-\sum_{\alpha=0}^{k-1}\bar u^1_\alpha u^2_\alpha+
\sum_{\alpha=k}^{k+l} \bar u^1_\alpha u^2_\alpha=-\langle\bfa^1,\bfa^2\rangle+
\langle\bfb^1,\bfb^2\rangle
\end{equation}
Here $\langle\bfa^1,\bfa^2\rangle$ denotes the standard quaternionic-Hermitian
inner product on $\HQ^{k}$ and we shall denote the associated norm by
$|\!|\bfa|\!|^2=\langle\bfa,\bfa\rangle$.  The form $F_{k,l}$ defines the flat
semi-Riemannian metric of signature $(4k,4l)$ on $\HQ^{k,l}$.

\begin{definition} Let $\cH_{k,l}(\epsilon)=\{(\bfa,\bfb)\in \HQ^{k,l}\ | \
-|\!|\bfa|\!|^2+|\!|\bfb|\!|^2=\epsilon\}$.
\begin{enumerate}
\item
$\cH_{k,l}(-1)\simeq S^{4k-1}\times\HQ^l$, where $k>0$, is a semi-Riemannian
submanifold of signature $(4k-1,4l)$ called the pseudosphere.
\item
$\cH_{k,l}(+1)\simeq \HQ^k\times S^{4l-1}$, where $l>0$, is a semi-Riemannian
submanifold of signature $(4k,4l-1)$ called the pseudohyperboloid.
\item
$\cH_{k,l}(0)\simeq S^{4k-1}\times S^{4l-1}\times\R/{\sim}$, where $k,l>0$ and
$\sim$ identifies $S^{4k-1}\times S^{4l-1}\times\{0\}$ with a point, is called
the null cone.
\end{enumerate}
\end{definition}
Let $k+l=n+1$ and $\Symp(k,l)\subset \GL(n+1,\HQ)$ which preserves the form
$F_{k,l}$. It is well-known that $\cH_{k,l}(\pm1)$ are spaces of constant
curvature and as homogeneous spaces of the semi-symplectic group $\Symp(k,l)$
they are
\begin{equation*}
\cH_{k,l}(\epsilon)=\begin{cases}
\Symp(k,l)/\Symp(k,l-1)& \text{when $l>0$ and $\epsilon=-1$},\\
\Symp(k,l)/\Symp(k-1,l)& \text{when $k>0$ and $\epsilon=+1$}. \end{cases}
\end{equation*}
Consider $\HQ_-^{k,l}=\{(\bfa,\bfb)\in\HQ^{k,l}\ \ |\ \ |\!|\bfa|\!|^2<
|\!|\bfb|\!|^2\},$ and $\HQ_+^{k,l}=\{(\bfa,\bfb)\in\HQ^{k,l}\ \ |\ \
|\!|\bfa|\!|^2> |\!|\bfb|\!|^2\}.$ Also, let us write $\HQ_0^{k,l}$ for
$\cH_{k,l}(0)$ as an alternative notation. We can then write
\begin{equation}\label{union}
\HQ^{k,l}=\HQ_-^{k,l}\cup \HQ_0^{k,l}\cup \HQ_+^{k,l}.
\end{equation}
After removing ${\bf 0}\in \HQ^{k,l}$ we consider the action of $\HQ^*$
on~\eqref{union} by right multiplication.

\begin{definition} Let $\HQ^{k,l}$ be the quaternionic vector space with
semi-hyperk\"ahler metric of signature $(4k,4l)$. We define the following
projective spaces.
\begin{enumerate}
\item $\HH{k,l-1}:=
\Proj_\HQ(\HQ_-^{k,l})=\HQ_-^{k,l}/\HQ^*=\cH_{k,l}(-1)/\Symp(1)$,
\item $\HH{k-1,l}:=
\Proj_\HQ(\HQ_+^{k,l})=\HQ_+^{k,l}/\HQ^*=\cH_{k,l}(+1)/\Symp(1)$,
\item $\Proj_\HQ(\HQ_0^{k,l})=(\HQ_0^{k,l}\setminus\{{\bf 0}\})/\HQ^*
=S^{4k-1}\times_{\Symp(1)}S^{4l-1}$.
\end{enumerate}
\end{definition}
If we make a choice of $\C^*\subset\HQ^*$ we also have complex
`projective' spaces.
\begin{definition} Let $\HQ^{k,l}$ be the quaternionic vector space with
semi-hyperk\"ahler metric of signature $(4k,4l)$. Let $\C^*\subset\HQ^*$.
We define
\begin{enumerate}
\item $\Proj_\C(\HQ_-^{k,l})=\HQ_-^{k,l}/\C^*=\cH_{k,l}(-1)/U(1)$,
\item $\Proj_\C(\HQ_+^{k,l})=\HQ_+^{k,l}/\C^*=\cH_{k,l}(+1)/U(1)$,
\item $\Proj_\C(\HQ_0^{k,l})=(\HQ_0^{k,l}\setminus\{{\bf 0}\})/\C^*
=S^{4k-1}\times_{U(1)}S^{4l-1}$.
\end{enumerate}
\end{definition}

\begin{proposition}
As homogeneous spaces of the semi-symplectic group
\begin{equation}\begin{split}
\Proj_\HQ(\HQ_-^{k,l})=\frac{\Symp(k,l)}{\Symp(1)\times \Symp(k-1,l)},&
\qquad \Proj_\C(\HQ_-^{k,l})=\frac{\Symp(k,l)}{U(1)\times \Symp(k-1,l)},
\ \ k>0\\
\Proj_\HQ(\HQ_+^{k,l})=\frac{\Symp(k,l)}{\Symp(k,l-1)\times \Symp(1)},&
\qquad
\Proj_\C(\HQ_+^{k,l})=\frac{\Symp(k,l)}{\Symp(k,l-1)\times U(1)},
\ \ l>0.
\end{split}\end{equation}
\end{proposition}

Furthermore, we have the natural fibrations
\begin{equation}
\begin{array}{ccccccc}&&\HQ_-^{k,l}&& \HQ_+^{k,l}&&\\
&\swarrow&&{}& &{\searrow}&\\ 
\Proj_\C(\HQ_-^{k,l}) && \decdnar{}&& \decdnar{}&&\Proj_\C(\HQ_+^{k,l}) \\
&\searrow& &{}& &{\swarrow}&\\
&&\Proj_\HQ(\HQ_-^{k,l})&& \Proj_\HQ(\HQ_+^{k,l})&&
\end{array}
\end{equation}
which can be glued together along  the common boundary
\begin{equation}
\begin{array}{ccccc}&\C^*&\rightarrow&\HQ_0^{k,l}\setminus\{{\bf 0}\}&\\
&&&\downarrow&\\
&S^2&\rightarrow&\Proj_\C(\HQ_0^{k,l})&\simeq S^{4k-1}\times_{S^1}S^{4l-1}\\
&&&\downarrow&\\
&&&\Proj_\HQ(\HQ_0^{k,l})&\simeq  S^{4k-1}\times_{S^3}S^{4l-1}
\end{array}
\end{equation}
to give $\HP {k+l-1}$, its twistor space $\CP {2k+2l-1}$ and the vector space
$\HQ^{k+l}\setminus\{{\bf 0}\}$.  Note that $\Proj_\HQ(\HQ_0^{k,l})\simeq
S^{4k-1}\times_{S^3}S^{4l-1}$ is both $S^{4k-1}$-bundle over $\HP {l-1}$ and
$S^{4l-1}$-bundle over $\HP {k-1}$. The following proposition is
straightforward.

\begin{proposition} The manifolds $\Proj_\HQ(\HQ_-^{k,l})$, $k>0$, are
semi-quaternion K\"ahler with holonomy group $\Symp(k-1,l)\cdot \Symp(1)$,
index $\nu=k-1$, negative scalar curvature, twistor space
$\Proj_\C(\HQ_-^{k,l})$, and Swann bundle $\HQ_-^{k,l}$; furthermore,
$\Proj_\HQ(\HQ_-^{k,l})$ is the quaternionic $\HQ^l$-bundle over the standard
quaternionic projective space $\HQ\Proj^{k-1}$ associated to the quaternionic
Hopf fibration.  The manifolds $\Proj(\HQ_+^{k,l})$, $l>0$, are
semi-quaternion K\"ahler with holonomy group $\Symp(k,l-1)\cdot \Symp(1)$,
index $\nu=k$, positive scalar curvature, twistor space
$\Proj_\C(\HQ_+^{k,l})$, and Swann bundle $\HQ_+^{k,l}$; furthermore,
$\Proj_\HQ(\HQ_+^{k,l})$ is the quaternionic $\HQ^k$-bundle over the standard
quaternionic projective space $\HQ\Proj^{l-1}$ associated to the quaternionic
Hopf fibration.  Topologically, $\HH{k,l-1}=\Proj_\HQ(\HQ_-^{k,l})$ and
$\HH{k-1,l}=\Proj_\HQ(\HQ_+^{k,l})$ are the components of
$\HQ\Proj^{k+l-1}\smallsetminus \Proj_\HQ(\HQ^{k,l}_0)$.
\end{proposition}
The bundle structure of $\Proj_\HQ(\HQ_-^{k,l})$ is the one associated to the
right quaternionic multiplication of the quaternionic vector space $\HQ^k$ by
the unit quaternions.  Explicitly, let $(\bfa,\bfb)\in \cH_{k,l}(-1)\subset
\HQ_-^{k,l}$.  Let us identify $\cH_{k,l}(-1)\simeq S^{4k-1}(1)\times\HQ^l$
via a map
\begin{equation*}
f(\bfa,\bfb)=(\bfv,\bfb)=\biggl( \frac{\bfa}{|\!|\bfb|\!|^2+1},\bfb\biggr).
\end{equation*}
Then $\sigma\in \Symp(1)$ acting on $S^{4k-1}(1)\times\HQ^l$ by
$(\bfv,\bfb)\rightarrow(\bfv\sigma,\bfb\sigma)$ gives the quotient which can
be identified with $\Proj_\HQ(\HQ_-^{k,l})$. Hence, $\Proj_\HQ(\HQ_-^{k,l})$
is an $\HQ^k$-bundle (quaternionic vector bundle) over $\HQ\Proj^{l-1}$
associated to the quaternionic Hopf bundle $S^3\rightarrow S^{4l-1}\rightarrow
\HQ\Proj^{l-1}$.

\begin{example}\label{ex1}
Let $(k,l)=(1,2)$. Then $\Proj_\HQ(\HQ_-^{1,2})$ is simply
the unit open 8-ball in $\HQ^2$. The boundary of this cell
$\Proj_\HQ(\HQ^{1,2}_0)= S^3\times_{S^3} S^7\simeq S^7$ is the unit sphere. The
space $\Proj_\HQ(\HQ_+^{1,2})$ is the $\HQ\simeq\R^4$ bundle over $\HP1\simeq
S^4$ associated to the quaternionic Hopf bundle $S^3\rightarrow S^7\rightarrow
S^4$.  Viewed another way $\Proj_\HQ(\HQ_+^{1,2})$ is a complement of the unit
8-ball in $\HQ^2$ with $\HP1\simeq S^4$ added in at infinity.
\end{example}

\begin{remark} Note that the map $\psi\colon\HQ^{k,l}\to \HQ^{l,k}$ defined
by
\begin{equation}
\psi(u_0,u_1,\ldots, u_{k-1}, u_k,\ldots, u_n)=(u_n,\ldots,u_k,u_{k-1},
\ldots,u_0)
\end{equation}
is the anti-isometry (or metric reversal) which induces anti-isometries
\begin{equation*}
\psi\colon\cH_{k,l}(\epsilon)\to \cH_{l,k}(-\epsilon),\text{ i.e., }
\psi\colon \HH{k,l}\to \HH{l,k}.
\end{equation*}
For example, $\Proj_\HQ(\HQ_-^{n+1,0})$ is diffeomorphic to $\HP n$ but has
negative-definite metric. It can be identified with $\Proj_\HQ(\HQ_+^{0,n+1})$
which is obviously the usual definition of $\HP n$ by changing the sign of the
metric.  As a result we can restrict our discussion only to the negative
scalar curvature spaces $\Proj_\HQ(\HQ_-^{k,l})$, $k>0$. This is not natural
if one talks about the projective space $\Proj_\HQ(\HQ_-^{n+1,0})$ but in this
paper we will mostly deal with the case $k<n+1$.
\end{remark}

We now describe the spaces $(\Proj_\HQ(\HQ_-^{k,l}), g^-_{k,l})$ in
inhomogeneous quaternionic coordinates. One needs $k$ quaternionic charts to
cover $\Proj_\HQ(\HQ_-^{k,l})$, namely
\begin{equation}
\cU_\beta=\{\bfu\in \Proj_\HQ(\HQ_-^{k,l})\ \ |\ \ u_\beta\not=0\},\ \ \
\beta=0,\ldots, k-1.
\end{equation}
On $\cU_\beta$ we write
\begin{equation}
\bfx^\beta=(x^\beta_1,\ldots,x^\beta_n)=(u_0u_\beta^{-1},\ldots,
u_{\beta-1}u_\beta^{-1},u_{\beta+1}u_\beta^{-1}
\ldots,u_nu_\beta^{-1})\in\HQ^n.
\end{equation}
Note that~\eqref{symform} implies that on $\cU_\beta$ we have
\begin{equation}
1-F_{k-1,l}(\bfx^\beta,\bfx^\beta)=1+\sum_{\alpha=1}^{k-1}|x^\beta_\alpha|^2-
\sum_{\alpha=k}^{n}|x^\beta_\alpha|^2=1/|u_\beta|^2>0.
\end{equation}
Let us denote $F_{k-1,l}$ simply by $\langle*,*\rangle_{k-1,l}$ with
the associated semi-norm $|\!|*|\!|_{k-1,l}$.  Then, on $\cU_\beta$,
$|\!|\bfx^\beta|\!|_{k-1,l}<1$ and the semi-quaternion K\"ahler metric
$g^-_{k,l}$ reads
\begin{equation}
g^-_{k,l}=\frac{1}{1-|\!|{\bf x^\beta}|\!|_{k-1,l}^2}
\biggl(|\!|d{\bf x^\beta}|\!|_{k-1,l}^2+
\frac{1}{1-|\!|{\bf x^\beta}|\!|_{k-1,l}^2}|\langle d{\bf x^\beta},
{\bf x^\beta}\rangle_{k-1,l}|^2\biggr).
\end{equation}
We will often refer to $\bfu=(u_0,u_1,\ldots u_n)$ as homogeneous coordinates
on $\Proj_\HQ(\HQ_-^{k,l})$.
\begin{example} It is clear that $\Proj_\HQ(\HQ_-^{1,n})=\HH n$ is simply
the unit ball in $\HQ^n$ with the quaternionic hyperbolic metric. In this case
$\cU_0$ is the only chart so we have global inhomogeneous coordinates
\begin{equation}
\bfx=(x_1,\ldots,x_n)=(u_1u_0^{-1},\ldots,u_nu_0^{-1})\in\HQ^n.
\end{equation}
with the positive definite hyperbolic metric
\begin{equation}
g=\frac{1}{1-|{\bf x}|^2}\biggl(|d{\bf x}|^2+
\frac{1}{1-|{\bf x}|^2}|\langle d{\bf x},{\bf x}\rangle|^2\biggr).
\end{equation}
\end{example}
These are not the only examples of semi-quaternion K\"ahler manifolds as we
shall see. However, many other examples can be obtained by taking
semi-quaternionic K\"ahler quotients of $\Proj_\HQ(\HQ_-^{k,l})$ by subgroups
of $\Symp(k,l)$. The quotient construction in the semi-Riemannian case works
in the similar way as in the Riemannian case. However, the zero-level set for
the moment map need not be a semi-Riemannian submanifold. For example, when
$G$ is a 1-parameter subgroup acting on a semi-quaternion K\"ahler manifold
$(M^{4n},g)$ of index $\nu$ then $N=\mu^{-1}(0)\subset M$ can have regions of
signature $(4\nu-3,4n-4\nu)$ and $(4\nu,4n-4\nu-3)$ separated by all points in
$M$ with $g(V,V)=0$, where $V$ is the vector field of the $G$-action on
$M$. Let us call these two regions by $N_-=\mu_-^{-1}(0)$ and
$N_+=\mu_+^{-1}(0)$. We have

\begin{theorem}
Let $(M^{4n},g)$ be a semi-Riemannian manifold with quaternionic index
$\nu$ and $G\subset {\rm Isom}_\Omega(M,g)$ be a one-parameter subgroup of
isometries of $M$ preserving the quaternion K\"ahler $4$-form $\Omega$. Let
$\mu:M\to \cV$ be the quaternion K\"ahler moment map for this action and let
$\mu_-^{-1}(0)\subset M$, $\mu_+^{-1}(0)\subset M$ be a semi-Riemannian
submanifolds of signature $(4\nu-3,4\nu)$ and $(4\nu,4\nu-3)$. If $G$ acts
freely and properly on $\mu_\pm^{-1}(0)$ the quotients $M_-=\mu_-^{-1}(0)/G$
and $M_+=\mu_+^{-1}(0)/G$ are semi-quaternion K\"ahler manifolds of dimension
$4n-4$ and quaternionic index $\nu-1$ and $\nu$, respectively.
\end{theorem}

The situation is even more complex when we choose an arbitrary $G\subset {\rm
Isom}_\Omega(M,g)$. In general, depending on how $G$ acts on $M$ one
should separate $\mu^{-1}(0)$ into submanifolds of signature
$(4\nu-3c,4n-4\nu-3d)$, where $c+d={\rm dim}(G)$ and one could expect
quotients of various quaternionic indices ranging from $0$ to ${\rm min}({\rm
dim}(G), \nu)$.

However, in this paper we shall focus our interest on the special case when
$M=\Proj_\HQ(\HQ_-^{k,l})$ or $M=\Proj_\HQ(\HQ_+^{k,l})$ with $k+l=3$ and ${\rm
dim}(G)=k+l-2=1$.  When $(k,l)=\{(0,3),(3,0)\}$ we are in the realm of the
$S^1$ reductions of $\HP2$, which have been already studied in
\cite{MR89m:53075}: the quotients are orbifold complex weighted projective
planes. In the case of $(k,l)=(1,2)$ we have two projective spaces one can
consider: $\Proj_\HQ(\HQ_-^{1,2})=\HH2$ and $\Proj_\HQ(\HQ_+^{1,2})=\HH{1,1}$.
However, as described in Example~\ref{ex1} these are two pieces of $\HP2$
cut along a 7-sphere. The choice of $G\subset \Symp(1,2)$, simultaneously
determines the quaternion K\"ahler reduction of both $\HH2$ and $\HH{1,1}$ by
$G$. In fact, the reduction depends only on the conjugacy classes of such
1-parameter subgroups in $\Symp(1,2)$. These, on the other hand, are given by
adjoint orbits in the Lie algebra $\symp(1,2)$. For each such adjoint orbit
$[\Delta]$ ($\Delta\in\symp(1,2)$) one can consider 1-parameter group
\begin{equation}
G(\Delta)=\{A\subset \Symp(1,2)\ \ |\ \ A=e^{\Delta t},\ \ \ t\in\R\}
\end{equation}
acting on $\HQ^{1,2}$ as a subgroup of $\Symp(1,2)\subset\GL(3,\HQ)$. This
action descends to an action on
\begin{equation}
\HP2= \HH2 \cup S^7\cup \HH{1,1}
\end{equation}
preserving the above decomposition and defining the semi-quaternion K\"ahler
moment maps.  Following Swann \cite{MR92c:53030}, it is convenient to consider
the semi-hyperk\"ahler moment map $\mu\colon\HQ^{1,2}\rightarrow {\Imag}(\HQ)$
and the corresponding decomposition of the Swann bundle.  We then write
\begin{equation}
\mu^{-1}_\Delta({\bf 0})=N_-(\Delta)\cup N_0(\Delta)\cup N_+(\Delta),
\end{equation}
where $N_\epsilon(\Delta)$ are restriction of $\mu_\Delta^{-1}({\bf 0})$ to
$\HQ^{1,2}_\epsilon$. As we shall see $N_-(\Delta)$ can be empty,
$N_+(\Delta)$ is never empty. Let $N_-(\Delta)$ be nonempty and suppose
\begin{equation*}
\Proj_\HQ(N_-(\Delta))=N_-(\Delta)/\HQ^*\subset \HH2
\end{equation*}
is a submanifold in the 8-ball $\HH2$. Further assuming that $G(\Delta)$ acts
freely and properly on $\Proj_\HQ(N_-(\Delta))$ we define the quotient
\begin{equation}
G(\Delta)\rightarrow \Proj_\HQ(N_-(\Delta))\rightarrow M_-(\Delta)=
G(\Delta)\backslash N_\epsilon(\Delta)/\HQ^*.
\end{equation}
It follows that the metric $g(\Delta)$ on $M_-(\Delta)$ obtained by inclusion
and submersion in the quotient construction is a complete SDE metric of
negative scalar curvature.  Its Swann bundle
$\cU(M_-(\Delta))=G(\Delta)\backslash N_-(\Delta)$ is a semi-hyperk\"ahler
manifold of index 1. Hence, for every $\Delta\in\symp(1,2)$ such that
$\Proj_\HQ(N_-(\Delta)) \subset \HH2$ and $G(\Delta)$ acts freely and properly
on it we get a negative SDE manifold $(M_-(\Delta),g(\Delta))$.  What remains
is to enumerate all possible adjoint orbits (this will be done in the next
section) and examine all the possible quotients (the following four sections).

The projectivisation $\Proj_\HQ(N_+(\Delta))\subset \Proj_\HQ(\HQ_+^{1,2})$,
in general, does not need to be a semi-Riemannian submanifold. Let
$V(\bfu)=\Delta\cdot\bfu$ be the vector field for the $G(\Delta)$-action on
$\Proj_\HQ(\HQ_+^{1,2})$.  Then the norm square of $V$ in the semi-Riemannian
metric $g^+_{1,2}$ can be negative, positive, or it can vanish.  Let
$\Proj^+_\HQ(N_+(\Delta))\subset \Proj_\HQ(\HQ_+^{1,2})$ be the subset on
which $g^+_{1,2}(V,V)>0$ while $\Proj^-_\HQ(N_+(\Delta))\subset
\Proj_\HQ(\HQ_+^{1,2})$ the subset on which $g^+_{1,2}(V,V)<0$.  If
$\Proj^+_\HQ(N_+(\Delta))$ is a submanifold in $\Proj_\HQ(\HQ_+^{1,2})$ then
it is a semi-Riemannian submanifold of signature $(4,1)$. On the other hand,
if $\Proj^-_\HQ(N_+(\Delta))$ is a submanifold in $\Proj_\HQ(\HQ_+^{1,2})$
than it is a semi-Riemannian submanifold of signature $(1,4)$. At least
locally we can define two different quotient metrics: (1) if
$\Proj^+_\HQ(N_+(\Delta))$ is not empty we have positive scalar curvature
metric $g^+(\Delta)$ on $\Proj^+_\HQ(N_+(\Delta))/G(\Delta)$ of signature
$(4,0)$ (anti-Riemannian); (2) if $\Proj^-_\HQ(N_+(\Delta))$ is not empty we
have positive scalar curvature metric $g^-(\Delta)$ on
$M_+^-(\Delta)=\Proj^-_\HQ(N_+(\Delta))/G(\Delta)$ of signature $(0,4)$. The
metric $g^-(\Delta)$, is a Riemannian metric of positive scalar
curvature. Typically this metric is not complete, unless the quotient can be
globally extended to the symmetric metric on $S^4$ or $\CP2$. On the other
hand $g^+(\Delta)$ is anti-Riemannian metric of positive scalar curvature so
that $-g^+(\Delta)$ is a Riemannian metric of negative scalar curvature on
$M_+^+(\Delta)=\Proj^+_\HQ(N_+(\Delta))/G(\Delta)$. Generally this metric is
not complete.  However, as we shall see in section~\ref{berg} complete metrics
of this type can occur.

Hence, a priori, for each $\Delta \in\symp(1,2)$ we have locally three
different metrics: $g(\Delta)$, $-g^+(\Delta)$ and $g^-(\Delta)$. The two
metrics $g(\Delta)$, $-g^+(\Delta)$ are negative SDE while $g^-(\Delta)$ is
positive SDE.

\begin{remark} Similarly, we can consider any orbits $[\Delta]$ under
$\Symp(1,n)$ of the $(n-1)$-dimensional subalgebras $\Delta\subset\symp(1,n)$.
Our analysis carried out for $\Symp(1,2)$ applies without any changes and, a
priori, for each $\Delta$ we obtain locally 3 different metrics: $g(\Delta)$,
$-g^+(\Delta)$ and $g^-(\Delta)$.  In addition, when $\Delta\subset\symp(1,n)$
is Abelian these metrics have two commuting Killing vectors.  Even more
generally, we could consider any orbit $[\Delta]$ under $\Symp(k,l)$ of
$(k+l-2)$-dimensional subalgebras $\g\subset\symp(k,l)$. If both $k,l$ are
greater than one $\Proj_\HQ(\HQ_\pm^{k,l})$ are both semi-quaternion K\"ahler.
Our analysis carried out for $(1,2)$ still applies and, a priori, for each
$\Delta$ we get locally four different metrics: two from the reduction of
$\Proj_\HQ(\HQ_-^{k,l})$ and the other two from the reduction of
$\Proj_\HQ(\HQ_+^{k,l})$. Two of these metrics will have negative scalar
curvature. The case of $k=l$ is of special interest as we shall see in
section~\ref{berg}. Again, for Abelian subalgebras the metrics will have two
commuting Killing vectors while the non-Abelian case is more general.
\end{remark}

\section{Adjoint orbits in $\symp(1,2)$}\label{sp(1,2)}

Adjoint orbits of elements in the classical Lie algebras $\g$ have been
determined by Burgoyne and Cushman \cite{MR55:5761}. We shall use this work to
find all the conjugacy classes of one-parameter subgroups of
$\Symp(1,2)$. First let us review some basic definitions. The symmetric form
on $\HQ^{1,2}$ is given by ${\mathbf u}^\dagger \F \mathbf u$, where
\begin{equation}
\mathbf u=\begin{pmatrix} u_0\\ u_1\\ u_2 \end{pmatrix}
\qquad\qquad
\F=\F_{1,2}=
\begin{pmatrix}
-1&0&0\\ 0&1&0\\ 0&0&1
\end{pmatrix}.
\end{equation}
We can describe $\Symp(1,2)$ and its Lie algebra $\symp(1,2)$ as $3\times3$
matrices preserving $\F$, i.e.,
\begin{align}
\Symp(1,2)&=\{g\in \cM_{3\times3}(\HQ)\ \ |\ \ g^\dagger\F g=\F\}\\
\symp(1,2)&=\{Y\in \cM_{3\times3}(\HQ)\ \ |\ \ \F Y+Y^\dagger\F={\bf 0}\}
\end{align}
Explicitly, an element of $\symp(1,2)$ can be written as
\begin{equation}
Y=\begin{pmatrix}
a &\alpha &\beta \\ 
\bar\alpha & b &\gamma \\ 
\bar\beta & -\bar\gamma& c \\
\end{pmatrix},
\end{equation}
Setting $\alpha=\beta=0$ gives the maximal compact subalgebra
$\symp(1)\oplus\symp(2)$ while $\beta=\gamma=0$ yields $\symp(1,1)\oplus
\symp(1)$.

We say that $Y\in\symp(1,2)$ is {\it decomposable} if $\HQ^{1,2}$ may be split
non-trivially as a direct sum of mutually orthogonal $Y$-invariant
quaternionic subspaces. Otherwise we say that $Y$ is {\it
indecomposable}. Choosing a particular unit quaternion $i$ identifies
$\HQ^{1,2}\cong \HQ^3$ with $\C^6$, which realizes $\symp(1,2)$ as a
subalgebra of $\gl(6,\C)$. We shall say an element $Y\in \symp(1,2)$ is {\it
semisimple} iff it is diagonalizable as an element of $\gl(6,\C)$. Any $Y\in
\symp(1,2)$ can be uniquely written as $Y=S+N$, where $S$ is semisimple, and
$N$ is nilpotent with $[S,N]=0$. If $N^{m+1}=0, N^{m}\neq0$ then the integer
$m$ is called the {\it height} of $Y$.  Semisimple elements have height equal
to zero.

\begin{definition}\label{def2}
We define the following elements of $\symp(1,2)$\textup:
\begin{align*}
T_0(ip_0,ip_1,ip_2)=&\ \  
\begin{pmatrix}
ip_0&0&0\\
0&ip_1&0\\
0&0&ip_2\\
\end{pmatrix},\\
T_0(\lambda,ip,iq)=&\ \
\begin{pmatrix}
ip &\lambda &0\\
\lambda &ip&0\\
0&0&iq\\
\end{pmatrix},\\
T_1(\lambda,ip,iq)=&\ \
\begin{pmatrix}
ip &0 &0\\
0 &ip&0\\
0&0&iq
\end{pmatrix}+
\lambda\begin{pmatrix}
i&i&0\\
-i&-i&0\\
0&0& 0
\end{pmatrix},\\
T_2(\lambda,ip)=&\ \ 
ip\,\bbi_3+
\lambda\begin{pmatrix}
0&0&-i\\
0&0&i\\
i&i&0
\end{pmatrix},
\end{align*}
where \textup(throughout\textup) $\lambda\neq 0$.
\end{definition}

The first two 3-parameter families of elements are semisimple and they are in
two different Cartan subalgebras of $\symp(1,2)$. They are necessarily
decomposable. The first one corresponds to the decomposition of $\HQ^{1,2}$
into $\HQ\oplus\HQ\oplus\HQ$ while the second decomposes $\HQ^{1,2}$ into
$\HQ^{1,1}\oplus\HQ$. The 3-parameter family $T_1(\lambda,ip,iq)$ has height
one (and $T_1(\lambda,0,0)$ is 2-step nilpotent). These are decomposable,
splitting $\HQ^{1,2}$ into $\HQ^{1,1}\oplus\HQ$.  Finally, the 2-parameter
family $T_2(\lambda,ip)$ has height two (and $T_2(\lambda,0)$ is $3$-step
nilpotent). These are indecomposable.  Note that all elements in the
Definition~\ref{def2} are inside the subalgebra $\un(1,2)$. Furthermore, note
that we chose $T_1:= T_1(1,0,0)$ and $T_2:= T_2(1,0)$, so that they
commute. In fact $\{i\bbi_3,T_2, T_1=iT_2^2\}$ span a maximal nilpotent
Abelian subalgebra of $\symp(1,2)$.

The following proposition follows from \cite{MR55:5761}.
\begin{proposition}
Let $Y$ be an arbitrary non-zero element of $\symp(1,2)$. Then $Y$ is
conjugate under the adjoint $\Symp(1,2)$ action to an element $\Delta$ of
Definition~\ref{def2}. This element is unique, except in the height one case,
where $T_1(\lambda,p,q)$ is conjugate to $T_1(1,p,q)$ or $T_1(-1,p,q)$ for
$p\neq 0$, and to $T_1(1,0,q)$ for $p=0$, and in the height two case, where
$T_2(\lambda,p)$ is conjugate to $T_2(1,p)$.

Furthermore, any one-parameter subgroup in $\Symp(1,2)$ is conjugate to
$G(\Delta)=\{A\in \Symp(1,2)\ \ | \ \ A=e^{\Delta t}\}$, where $\Delta$ is one
of the types of the Definition~\ref{def2}.
\end{proposition}
In other words, the list of Definition~\ref{def2} enumerates all adjoint
orbits in $\symp(1,2)$. The corresponding conjugacy classes of one parameter
subgroups of $\Symp(1,2)$ are enumerated by these elements up to scale:
$\Delta$ and $c\Delta$ define the same subgroup for any $c\neq 0$.

In the following, it will sometimes be more convenient to work with a
different basis of $\HQ^{1,2}$ in which the symmetric form may be
written ${\mathbf v}^\dagger \tilde\F \mathbf v$ with
\begin{equation}
\mathbf v=\begin{pmatrix} v_0\\ v_1\\ v_2 \end{pmatrix}
\qquad\qquad
\tilde\F=
\begin{pmatrix}
0&1&0\\ 1&0&0\\ 0&0&1
\end{pmatrix}.
\end{equation}
The advantage of this basis is that the last three matrices in
Definition~\ref{def2} are conjugated to the following simpler forms.
\begin{align*}
\tilde T_0(\lambda,ip,iq)=&\ \
\begin{pmatrix}
ip +\lambda & 0 &0\\
0 &ip-\lambda&0\\
0&0&iq
\end{pmatrix},\\
\tilde T_1(\lambda,ip,iq)=&\ \
\begin{pmatrix}
ip &0 &0\\
-i\lambda &ip&0\\
0&0&iq
\end{pmatrix},\\
\tilde T_2(\lambda,ip)=&\ \ 
ip\,\bbi_3+
\lambda\begin{pmatrix}
0&0&0\\
0&0&i\\
i&0&0
\end{pmatrix}.
\end{align*}

We end our discussion by noting that it is straightforward to compute the
momentum map in homogeneous coordinates associated to a generator $T$ or
$\tilde T$ using the general formulae
\begin{align*}
\mu_T(\bfu) &= \bfu^\dag \F T \bfu = (-\bar u_0, \bar u_1, \bar u_2) T
\Biggl(\begin{matrix} u_0\\ u_1\\ u_2 \end{matrix}\Biggr)\\
\mu_{\tilde T}(\bfv) &= \bfv^\dag \tilde \F \tilde T \bfv =
(\bar v_1, \bar v_0, \bar v_2) \tilde T
\Biggl(\begin{matrix} v_0\\ v_1\\ v_2 \end{matrix}\Biggr).
\end{align*}

\section{The Pedersen--LeBrun Metrics on Line Bundles over $\CP1$}

In this section we will examine the case of $\Delta=\Delta_0(\bfp)=
T_0(ip_0,ip_1,ip_2)$. We shall assume that this generates a circle action,
which means, after rescaling $\Delta$, that we may assume that the $p_i$'s are
integers with ${\rm gcd}(p_0,p_1,p_2)=1$. (We can assume that the weights do
not vanish as the cases when one or two of the weights vanish are degenerate.)
We then have a circle action on the quaternionic hyperbolic 2-ball $\HH2$
given in homogeneous coordinates by
\begin{equation}\label{3.1}
\varphi_t(u_0,u_1,u_2)=(e^{2\pi ip_0t}u_0,\ e^{2\pi ip_1t}u_1,
\ e^{2\pi ip_2t}u_2)
\end{equation}
where $t\in [0,1)$.  We note that this action is effective unless the weights
$p_0,p_1,p_2$ are all odd, in which case we obtain an effective action of a
quotient circle by taking $t\in [0,1/2)$.  In inhomogeneous coordinates
$(x_1,x_2)$ we have
\begin{equation}
\varphi^\bfp_t(x_1,x_2)=(e^{2\pi ip_1t}x_1 e^{-2\pi ip_0t},
e^{2\pi ip_2t}x_2 e^{-2\pi ip_0t}),
\end{equation}
and the moment map is given as
\begin{align}
\mu_\bfp(\bfu)&=-p_0\bar u_0iu_0+p_1\bar u_1iu_1+p_2\bar u_2iu_2,\\
f_\bfp(\bfx)&=u_0\mu_\bfp(\bfu)u_0^{-1}=
-ip_0+p_1\bar x_1ix_1+p_2\bar x_2ix_2
\end{align}
in homogeneous or inhomogeneous coordinates. We now write
\begin{equation}\bfx=\bfz+\bfw j = \bfz+j\bar \bfw\end{equation}
where $\bfz,\bfw\in \C^2$  and observe that
\begin{equation*}
\varphi^\bfp_t\begin{pmatrix}z_1&w_1\\ z_2&w_2\end{pmatrix}=
\begin{pmatrix}e^{2\pi i(p_1-p_0)t}z_1&e^{2\pi i(p_1+p_0)t}w_1\\
e^{2\pi i(p_2-p_0)t}z_2&e^{2\pi i(p_2+p_0)t}w_2\end{pmatrix}.
\end{equation*}
\begin{equation}
\mu_\bfp^{-1}(0)=\lbrace (\bfz,\bfw)\in \HH2:
\ \sum_{\alpha=1,2}p_\alpha(|z_\alpha|^2-|w_\alpha|^2)=p_0,
\ \sum_{\alpha=1,2}p_\alpha\bar w_\alpha\cdot z_\alpha= 0\rbrace.
\end{equation}

\begin{proposition}\label{prop3} Let $q_\alpha=p_\alpha/p_0$. Then the subset
$\mu_\bfp^{-1}(0)\subset \HH2$ is empty unless $|q_\alpha|>1$ for at least one
$\alpha$. Otherwise, $\mu_\bfp^{-1}(0)$ is an open smooth submanifold of
codimension $3$.
\end{proposition}

\begin{proof} We can assume that both $q_1,q_2$ are positive (otherwise
we simply reverse the role of $z_\alpha$ and $w_\alpha$ in the argument
below).  On the one hand, the momentum constraint gives
\begin{equation*}
q_1(|z_1|^2-|w_1|^2)+q_2(|z_2|^2-|w_2|^2)=1,
\end{equation*}
so that
\begin{equation*}
-q_1|z_1|^2-q_2|z_2|^2\leq -1.
\end{equation*}
On the other hand, we have
\begin{equation*}
|z_1|^2+|z_2|^2\leq |z_1|^2+|w_1|^2+|z_2|^2+|w_2|^2<1
\end{equation*}
by the unit ball condition. Adding the two inequalities we get
\begin{equation*}
(1-q_1)|z_1|^2+(1-q_2)|z_2|^2<0.
\end{equation*}
This has no solutions when $1\geq q_1$ and $1\geq q_2$. Otherwise, if (say)
$|q_1|>1$ then by taking $z_2,w_2=0$, it is easy to see that
$\mu_\bfp^{-1}(0)$ is nonempty. The last statement follows because
straightforward computation reveals that $0$ is a regular value of the
Jacobian of $\mu_\bfp$.
\end{proof}

Without loss of generality we will further assume that all weights are
positive.  We will also choose $q_1>1$ that is that $p_1>p_0$. Then we
have the following

\begin{proposition}\label{p:freeness} The $\varphi_t^\bfp$-action on the
level set $\mu_\bfp^{-1}(0)\subset \HH2$ is free if and only if $p_1=p_0+1$
and $0<p_2\leq p_0+1$ when one of the weights is even, or $p_1=p_0+2$ and
$0<p_2\leq p_0+2$ when all the weights are odd.
\end{proposition}

\begin{proof} Consider the set described by $(z_1,0,0,0)$. This meets
$\mu_\bfp^{-1}(0)$ in a circle, but any point on this circle is fixed by
$\Z_{p_1-p_0}$. Hence we must have $p_1=p_0+1$, unless all weights are odd
when we have $p_1=p_0+2$. Next, suppose that $p_2>p_0$. Then the set described
by $(0,z_2,0,0)$ also meets $\mu_\bfp^{-1}(0)$ in a circle and any point on
this circle is fixed by $\Z_{p_2-p_0}$. Thus if $p_2>p_0$, we must have $p_2 =
p_0+1$ (or $p_0+2$ if all weights are odd).  It is easy to see that $p_2$ can
be any integer with $0<p_2\leq p_0+1$ (or $p_0+2$ if all weights are odd).
\end{proof}

We now have:

\begin{theorem}\label{PLeB} For $\bfp\in\Z_+^3$ as in
Proposition~\textup{\ref{p:freeness}}, the quotient
$M(\bfp)=\mu_\bfp^{-1}(0)/S^1(\bfp)$ is a complete self-dual Einstein manifold
with negative scalar curvature and at least two commuting Killing
vectors. When $p_2=p_1$ \textup(which is $p_0+1$ or $p_0+2$\textup) the metric
is $U(2)$-invariant while when $p_2=p_0$ the metric is $U(1,1)$-invariant.
\end{theorem}

\begin{proof}
Only completeness of the induced metric on $M(\bfp)$ remains to be proven, and
this follows from the fact that the induced metric on the closed embedded
submanifold $\mu^{-1}_\bfp(0)\hookrightarrow \HH2$ is complete and the action
\ref{3.1} is proper.
\end{proof}

We continue with describing the total space of these metrics.  When $p_2=p_1$,
we expect that the metric is complete and, hence, it has to be one of the
possibilities listed by Hitchin in Theorem 13 of \cite{MR96g:53057}.  We will
show that our quotient metrics are the Pedersen--LeBrun metrics on complex line
bundles $\cO(n)\to \CP1$, $n\geq3$ (Theorem 13:3(d) of \cite{MR96g:53057}).
Before we analyze $M(p,p+1,p+1)$ and $M(p,p+2,p+2)$ let us recall a standard
description of a complex line bundle over $\CP1$ with first Chern class
$s$. Let $S^{3}=\lbrace \bfv\in \C^2:\ \ \|\bfv\|=1\rbrace$ and let
$s\in\Z^+$.  Then we set
\begin{equation}
{\cL}_{s}\equiv S^3\times \C/\Phi^{s},
\end{equation}
where $\Phi^{s}$ is the action of $S^1$ on $S^3\times \C$ given by
\begin{equation}
\Phi^{s}_\tau(\bfv,\alpha)=(\tau\bfv,\ \tau^s\alpha).
\end{equation}
The natural projection ${\cL}_{s}\ \longrightarrow S^2\cong S^3/S^1$ makes
${\cL}_{s}$ a complex line bundle over $S^2$ with $c_1({\cL}_{s})=s$.

Note that we get the same conclusion when we replace $\C$ by
$D_\C^1(1)=\{\alpha\in\C :|\alpha|<1\}.$ Then ${\cL}_{r,s}$ is a complex unit
disk bundle with first Chern class $s$.  Now we are ready for

\begin{theorem} Let $p\in\Z$ and $\bfp=(p-1,p,p)$, $p>1$. Then the quotient
metric $g(\bfp)$ is complete, $U(2)$-invariant and the total space $M(\bfp)$
can be identified with the complex unit disk bundle $\cL_{2p}\rightarrow\CP1$
with first Chern class equal to $2p$.
\end{theorem}

\begin{proof}
By Theorem \ref{PLeB} it suffices only to identify the quotient in this
special case.  Let $(\bfz,\bfw)\in{\HH2}$.  We make a slight change of these
coordinates by setting
\begin{equation}\bfx=
\frac{1}{\sqrt{p_0/p+\|\bfw\|^2}}\bfz,\ \ \ \ \bfy=\sqrt{2p}\bfw.
\end{equation}
In these coordinates the moment map equations can be written
\begin{equation}\mu^{-1}_\bfp(0)=\lbrace (\bfx, \bfy)\in \C^2\times \C^2:\ \
\ \ \|\bfx\|^2= 1,\ \ \ \bar \bfy\cdot \bfx=0,
\ \ ||\bfy||<1\rbrace.
\end{equation}
Then the circle action is given by
\begin{equation}\label{3.13}
\varphi_\tau(\bfx, \bfy)= (\tau\bfx,\ \tau^{2p-1}\bfy)
\end{equation}
for $\tau=e^{2\pi it}\in S^1$. Consequently, we have that $M(\bfp)$ is
equivalent to the quotient of the set
\begin{equation*}
\mu_\bfp^{-1}(0)=\lbrace(\bfx,\bfy)\in S^3\times D_\C^2(1):
\ \ \bfx\bot \bfy\rbrace\simeq S^3\times D_\C^1(1)
\end{equation*}
by the action~\eqref{3.13}. We define a map $f:S^3\times D_\C^1(1)\
\longrightarrow M(\bfp)$ by setting $f(\bfv,\alpha)=(\bfv,\alpha
\bfv^\dagger)$ where if $\bfv=(v_o,v_1)$ then $\bfv^\dagger=(-\bar v_1,
\bar v_o)$.  Note that for any $\tau \in S^1$ we have the commutative
diagram
\begin{equation*}
\begin{matrix}
(\bfv,\alpha)&&\longrightarrow&&(\bfv,\alpha \bfv^\dagger)\\
 \Big\downarrow &&&&\Big\downarrow\\
(\tau\bfv,\tau^{2p}\alpha)&&\longrightarrow&&(\tau\bfv,\tau^{2p-1}\alpha
\bfv^\dagger)
\end{matrix};
\end{equation*}
i.e., we have that $f\circ \Phi_\tau^{2p}=\varphi_\tau\circ f$. Thus $f$ is an
$S^1$-equivariant diffeomorphism and therefore $f$ induces a smooth
equivalence of the quotient spaces. Hence $M(\bfp)\simeq \cL_{2p}$.

\end{proof}

When $\bfp=(p-2,p,p)$ we immediately get the other half of the line bundles
with odd Chern classes:

\begin{theorem} Let $p\in\Z$ and $\bfp=(p-2,p,p)$, $p=2k+1>2$.  Then the
quotient metric $g(\bfp)$ is complete, $U(2)$-invariant and the total space
$M(\bfp)$ can be identified with the complex unit disk bundle
$\cL_{p}\rightarrow\CP1$ with first Chern class equal to $p$.
\end{theorem}

Note that this construction does not give the line bundles over $\CP1$ with
first Chern classes $c_1=1,2$. The metrics on $\cL_p$ with $p\geq 3$ have a
curious history. The quotient construction presented here was written in
\cite{MR88j:53076}.  The Pedersen metric on the $4$-ball~\cite{MR87i:53070}
depends on a single parameter $m^2\in(-1,\infty)$. It was realized later (see
\cite{MR96g:53057}) that setting this parameter to $(2-n)/n$ (with
$n\in\Z,n>2$) allows for the analytic continuation of this metric to a
complete metric on $\cO(n)\to \CP1$. The reason these metrics are called
Pedersen--LeBrun in \cite{MR96g:53057} is that they are conformal to the scalar
flat K\"ahler metrics on $\cO(-n)\to\CP1$ constructed by LeBrun
\cite{MR89f:53107}.

When, $p_0+1=p_1>p_2>0$ we take a different approach.  Observe that one can
still easily solve the complex equation of the moment map by setting
\begin{equation}
(w_1,w_2)=\alpha(-p_2\bar z_2, p_1\bar z_1),
\end{equation}
where $\alpha\in\C$.  The unit ball condition in terms of
$(z_1,z_2,\alpha)$ reads:
\begin{equation}|z_1|^2+|z_2|^2+|\alpha|^2p_1^2|z_1|^2+|\alpha|^2p_2^2|z_2|^2=
|z_1|^2(1+p_1^2|\alpha|^2)+|z_2|^2(1+p_2^2|\alpha|^2)<1,\end{equation}
while the remaining moment map equation is
\begin{equation*}
|z_1|^2(p_1-p_2p_1^2|\alpha|^2)+|z_2|^2(p_2-p_1p_2^2|\alpha|^2)=p_0.
\end{equation*}
Let us solve this equation with respect to $|z_1|^2$:
\begin{equation}|z_1|^2=\frac{p_0}{p_1}
\frac{1}{1-p_1p_2|\alpha|^2}-\frac{p_2}{p_1}|z_2|^2.
\end{equation}
One can immediately see that $z_1$ cannot vanish as then
\begin{equation}
|z_2|^2=\frac{p_0}{p_2}\frac{1}{1-p_1p_2|\alpha|^2}\geq\frac{p_0}{p_2}\geq1.
\end{equation}
Let $\rho=\frac{z_1}{|z_1|}$. It is easy to see that
\begin{equation*}
\phi_\tau(\rho,z_2,\alpha)=(\tau\rho,\tau^{p_2-p_0}z_2,\tau^{p_2+p_1}\alpha).
\end{equation*}

\begin{proposition} The level set $\mu^{-1}_\bfp(0)\simeq D\times S^1$,
where $D\subset\C^2$ is an open $4$-ball.
\end{proposition}

\begin{proof}
It is clear that $(\rho,z_2,\alpha)\in S^1\times D$ are coordinates on
$\mu^{-1}_\bfp(0)$. We have to check that $D$ is diffeomorphic to a 4-ball. To
do that let us consider
\begin{equation*}
\biggl(\frac{p_0}{p_1}\frac{1}{1-p_1p_2|\alpha|^2}
-\frac{p_2}{p_1}|z_2|^2\biggr)
(1+p_1^2|\alpha|^2)+|z_2|^2(1+p_2^2|\alpha|^2)<1
\end{equation*}
which can be written as
\begin{equation*}
f_{\bfp}(z_2,\alpha)=
(p_1-p_2)|z_2|^2[1-p_1p_2|\alpha|^2]^2+p_1^2p_2|\alpha|^2-1<0.
\end{equation*}
One can easily see that $|\alpha|^2<\frac{1}{p_1^2p_2}$
\begin{equation*}
D(\bfp)=\{(z_2,\alpha)\in\C\times\C \ \ | \ \ f_{\bfp}(z_2,\alpha)<0\}
\end{equation*}
is an open 4-ball.
\end{proof}

\begin{theorem} The quotient $M(\bfp)\simeq D(\bfp)$ is diffeomorphic to
a $4$-ball. The self-dual Einstein metric $g(\bfp)$ obtained from the
quaternion K\"ahler quotient is complete and it has two commuting Killing
vectors. Furthermore, $M(p,p+1,p)$ is of cohomogeneity one with respect to
$U(1,1)$.
\end{theorem}

The cohomogeneity one $U(1,1)$ action on $M(p,p+1,p)$ can be explicitly
described as follows. Let
\begin{equation*}
\bba=
\begin{pmatrix}
a&\tau b\\
\bar{b}&\tau\bar{a}\\
\end{pmatrix}\in U(1,1),
\end{equation*}
where $a,b,\tau\in\C$ with $|a|^2-|b|^2=1$ and $|\tau|^2=1$. This group acts
on the quaternionic ball as
\begin{equation*}
\varphi_\bba(\bfu)=
\begin{pmatrix}
a&0&\tau b\\
0&1&0\\
\bar{b}&0&\tau \bar{a}\\
\end{pmatrix}
\begin{pmatrix}
u_0\\ u_1\\ u_2\\
\end{pmatrix}
\end{equation*}
and, it commutes with the circle action given by $\varphi^\bfp_t$. In the
inhomogeneous chart we get
\begin{equation*}
\varphi_\bba
\begin{pmatrix}
x_1\\ x_2\\
\end{pmatrix}=
\begin{pmatrix}
x_1(a+\tau bx_2)^{-1}\\
(\bar{b}+\tau\bar{a}x_2)(a+\tau bx_2)^{-1}\\
\end{pmatrix}.
\end{equation*}
The above action preserves the zero level set of the moment map and it
descends to a cohomogeneity one isometric action on the quotient space
$M(p,p+1,p)$.  Cohomogeneity one SDE metrics with an isometric action of a
four-dimensional Lie group have been studied by Derdzi\'nski. Hence, up to
isometries $M(p,p+1,p)$ must be the cohomogeneity one self-dual K\"ahler
metric introduced in~\cite{Der:em} and more recently studied by Apostolov and
Gauduchon in \cite{math.DG/0003162}.

\section{Generalized Pedersen Metrics on the Ball}

In this section and the following two, we will consider the $\R$-actions on
$\HH2$ whose generators do not belong to the Lie algebra of a maximal
torus. To do this we shall work in the $\bfv=(v_0,v_1,v_2)$ coordinates
introduced in section~\ref{sp(1,2)}. In these coordinates $\HH2$ is the open
subset of $\HP2$ defined by the equation
\begin{equation}
{\bar v_0} v_1 + {\bar v_1}v_0 + |v_2|^2<0.
\end{equation}
It follows that $v_0$ does not vanish on $\HH2$ and so the inhomogeneous
coordinates $y_1=v_1v_0^{-1}$, $y_2=v_2 v_0^{-1}$ provide a global chart
identifying $\HH2$ with the domain
\begin{equation}
y_1+\bar y_1 + |y_2|^2 <0
\end{equation}
in $\HQ^2$. We remark (for later use) that the real part of $y_1$ is strictly
negative on this domain.

We begin by considering the case of $\Delta_0(p,q)=\tilde T_0(1,ip,iq)$, where
we have taken $\lambda=1$ by rescaling. The $\R$-action on the quaternionic
hyperbolic 2-ball $(\HH2, g)$ is given explicitly by
\begin{equation}
\varphi_t^{p,q}(\bfv)=
\begin{pmatrix}
e^{(ip+1)t}&0&0\\ 
0&e^{(ip-1)t}&0\\
0& 0 & e^{iq t}
\end{pmatrix}
\begin{pmatrix}v_0\\ v_1\\ v_2\end{pmatrix}
=\begin{pmatrix}e^{ipt}e^t v_0\\ e^{ipt}e^{-t}v_1\\ e^{iqt}v_2\end{pmatrix},
\end{equation}
which reduces, in inhomogeneous coordinates $\bfy=(y_1,y_2)$, to
\begin{equation}\label{eq:ped-act}
\varphi_t^{p,q}\begin{pmatrix}y_1\\ y_2\end{pmatrix}=
\begin{pmatrix}
e^{ip t} e^{-2t} y_1 e^{-ip t}\\
e^{iq t} y_2 e^{-ipt}
\end{pmatrix}.
\end{equation}

This action is a quaternionic isometry of the hyperbolic metric $g$ and it
defines a bundle valued momentum map $\mu_{p,q}\colon\HH2\rightarrow\cV$ given
in homogeneous and inhomogeneous coordinates by the function
\begin{align*}
\mu_{p,q}(\bfv) &= {\bar v_1}v_0-\bar v_0 v_1
+ p({\bar v_1} i v_0 + {\bar v_0} i v_1) + q{\bar v_2} i v_2,\\
f_{p,q}(\bfy) &= {\bar y_1}-y_1 + p({\bar y_1} i + i y_1) + q{\bar y_2} i y_2.
\end{align*}
Although this function is not invariant under the action~\eqref{eq:ped-act},
its zero set is, and the quaternion K\"ahler reduction of $\HH2$ by the one
parameter group $e^{\Delta_0(p,q) t}$ is the quotient of this zero set by the
group action. The resulting SDE metrics were first introduced in
\cite{MR88j:53076} and may be regarded as a deformation of the Pedersen
metrics on the ball to metrics with fewer symmetries.

\begin{theorem} Let $\Delta=\Delta_0(p,q)=\tilde T_0(1,ip,iq)$
and consider the one parameter group $\varphi_t^{p,q}=e^{\Delta t}$ acting on
the quaternionic hyperbolic space $\HH2$.  The quaternion K\"ahler reduction
$M(p,q)=\mu^{-1}_{p,q}(0)/\varphi^{p,q}$ is diffeomorphic to an open $4$-ball
for all $(p,q)\in\R^2$.  The quotient metric is complete, self-dual, and
Einstein of negative scalar curvature whose isometry group contains a
$2$-torus. Furthermore, the quotient metrics on $M(0,q)$ are isometric to the
Pedersen metrics, and their isometry group contains $U(2)$.
\end{theorem}

\begin{proof} Consider the following set
\begin{equation}
\cS_{p,q}=\{\bfy\ \ |\ \ f_{p,q}(y_1,y_2)=0 \text{ and }
y_1+{\bar y_1}=2\Real(y_1)=-1\}.
\end{equation}
For any $y_2\in\HQ$ there is a unique $y_1$ such that $\bfy\in\cS_{p,q}$.  It
follows that $\cS_{p,q}\cap\HH2$ is diffeomorphic to the open $4$-ball
$|y_2|^2 <1$ in $\HQ\cong \R^4$. Furthermore, $\cS_{p,q}\cap\HH2$ provides a
global slice for the action of $e^{\Delta t}$ in the zero set of the momentum
map: to see this, we only have to note that $e^{\Delta t}$ sends $y_1+{\bar
y_1}$ to $e^{-2t}(y_1+{\bar y_1})$ and therefore, since $y_1+{\bar y_1}<0$,
there is a unique $(y_1,y_2)$ in any orbit with $y_1+{\bar y_1}=-1$.

Therefore $M(p,q)$ is diffeomorphic to the open $4$-ball equipped with the
metric obtained by restriction to $\mu_{p,q}^{-1}(0)$ and submersion. The fact
that the quotient is a complete Riemannian manifold follows as in the proof of
Theorem \ref{PLeB}.  Moreover, it must be an SDE metric of negative scalar
curvature since it is obtained as quaternion K\"ahler quotient of $\HH2$.

The isometry group contains a $2$-torus since $\cS_{p,q}$ is invariant
under the transformation
\begin{equation*}
(y_1,y_2)\mapsto (\sigma y_1\sigma^{-1}, \tau
y_2\sigma^{-1})
\end{equation*}
by $(\tau,\sigma)\in U(1)\times U(1)$. As this action is by quaternionic
isometries on $\HH2$ and commutes with $e^{\Delta t}$, it descends to give an
action by isometries on $M(p,q)$. If $p=0$ this action may be extended, by
taking $\sigma\in Sp(1)$, to yield an action of $U(1)\cdot \Symp(1)\simeq
U(2)$.

To identify $M(0,q)$ as the Pedersen family, one can compute the metric
explicitly. Alternatively, we can use the classification of SDE metrics with
$\SU(2)$ symmetry by Hitchin \cite{MR96g:53057}. This classification provides
very few possible candidates with $U(2)$ symmetry: apart from the real and
complex hyperbolic metrics, the Pedersen metrics are the only examples. In
fact, one can see that $M(0,0)$ is real hyperbolic space but for other values
of $q$ the metric is not symmetric.
\end{proof}

\section{The Height One Quotients}

In this section we will examine the family of quotients of $\HH2$, obtained
from $\tilde T_1(\lambda,ip,iq)$. By rescaling we can assume that $p$ is $0$
or $1$, and if $p=0$ we can scale $q$ to $1$ or $0$. Since we assume $\lambda$
is nonzero, we can then conjugate so that $\lambda=\pm1$ (or $\lambda=1$ if
$p=0$) and rescale by the sign. Hence we only need to consider the quotients
$\Delta_1(p,q)=\tilde T_1(1,ip,iq)$ with $p\in\{-1,0,1\}$, and if $p=0$ we can
suppose $q\in\{0,1\}$. Nevertheless, for convenience we shall carry out our
analysis for arbitrary $p,q$. We have
\begin{equation}
\varphi_t^{p,q}(\bfv)=
\begin{pmatrix}
e^{ipt}&0&0\\ 
-it&e^{ipt}&0\\
0& 0 & e^{iq t}
\end{pmatrix}
\begin{pmatrix}v_0\\ v_1\\ v_2\end{pmatrix}
=\begin{pmatrix}e^{ipt}v_0\\ e^{ipt}v_1-it v_0 \\ e^{iqt}v_2\end{pmatrix},
\end{equation}
which reduces, in inhomogeneous coordinates $\bfy=(y_1,y_2)$, to
\begin{equation}\label{eq:h1-act}
\varphi_t^{p,q}\begin{pmatrix}y_1\\ y_2\end{pmatrix}=
\begin{pmatrix}
e^{ip t} (y_1-it) e^{-ip t}\\
e^{iq t} y_2 e^{-ipt}
\end{pmatrix}.
\end{equation}
The moment map for this action is given in homogeneous or inhomogeneous
coordinates by
\begin{align*}
\mu_{p,q}(\bfv) &= -\bar v_0 i v_0 + p({\bar v_0} i v_1+{\bar v_1} i v_0)
+ q{\bar v_2} i v_2,\\
f_{p,q}(\bfy)&= -i + p(i y_1+{\bar y_1} i) + q{\bar y_2} i y_2.
\end{align*}

\begin{theorem} Let $\Delta =\Delta_1(p,q)=\tilde T_1(1,ip,iq)$ and consider
the one parameter group $\varphi^{p,q}_t=e^{\Delta t}$ acting on the
quaternionic hyperbolic space $\HH2$. Then
\begin{enumerate}
\item
the quaternion K\"ahler reduction $M(p,q)=\mu^{-1}_{p,q}(0)/\varphi^{p,q}$ is
diffeomorphic to $\R^4$ for all $(p,q)$ with $p<0$.
\item 
the quaternion K\"ahler reduction $M(p,q)=\mu^{-1}_{p,q}(0)/\varphi^{p,q}$ is
diffeomorphic to $S^1\times \R^3$ for all $(p,q)$ with
$0\leq p<|q|$.
\end{enumerate}
In these cases $M(p,q)$ has a complete self-dual Einstein metric of negative
scalar curvature and its isometry group contains a $2$-torus.  In all other
cases \textup(i.e., if $p\geq|q|$\textup) the zero set of the momentum map is
empty.
\end{theorem}

\begin{proof} We begin by defining the set
\begin{equation}
\cS_{p,q}=\{(y_1,y_2)\ \ |\ \ f_{p,q}(y_1,y_2)=0\text{ and }
iy_1-{\bar y_1} i = 2\Real(iy_1)=0\}
\end{equation}
and claim that $\cS_{p,q}\cap \HH2$ can be identified with the quotient space
$M(p,q)$ as a global slice for the $\varphi^{p,q}$ action on the momentum zero
set.  Indeed, it is clear that as the action of $e^{\Delta t}$ sends
$\Real(iy_1)$ to $\Real(iy_1)+2t$, so there is a unique point of $\cS_{p,q}$ on
each orbit of $e^{\Delta t}$ in $\mu_{p,q}^{-1}(0)$. It remains to describe
the set $\cS_{p,q}\cap\HH2$.

For $p\neq 0$, there is a unique $(y_1,y_2)\in \cS_{p,q}$ for any
$y_2\in\HQ$. We now note that $\HH2$ is the domain
\begin{align*}
p(y_1 + {\bar y_1}) + p{\bar y_2} y_2 &< 0,\qquad p>0\\
p(y_1 + {\bar y_1}) + p{\bar y_2} y_2 &> 0,\qquad p<0.
\end{align*}
On the other hand
\begin{equation*}
0=\Real (i f_{p,q})= 1 - p (y_1+{\bar y_1}) - q \Real(i{\bar y_2} i y_2)
\end{equation*}
so that $\cS_{p,q}\cap\HH2$ may be identified with the set of $y_2\in \HQ$
satisfying
\begin{align*}
-q \Real(i{\bar y_2} i y_2) + p{\bar y_2} y_2 &< -1,\qquad p>0\\
-q \Real(i{\bar y_2} i y_2) + p{\bar y_2} y_2 &> -1,\qquad p<0.
\end{align*}
Writing $y_2=z_2+j w_2$ for $w_2,z_2\in \C$, this is the domain in $\C^2$
given by
\begin{align*}
(p+q)|z_2|^2+(p-q)|w_2|^2 &< -1,\qquad p>0\\
(p+q)|z_2|^2+(p-q)|w_2|^2 &> -1,\qquad p<0.
\end{align*}
For $p>0$, this domain is empty unless $p<|q|$, in which case it is the
exterior of a hyperboloid, which is diffeomorphic to $S^1\times\R^3$.  For
$p<0$, this domain is the interior of a hyperboloid for $-|q|<p<0$, the
interior of a cylinder for $p=-|q|$, and the interior of an ellipsoid for
$p<-|q|$: all these domains are diffeomorphic to $\R^4$.

We now consider the case $p=0$, when $y_1$ is not uniquely determined by
$y_2$. For $q=0$, the momentum zero set is empty. Otherwise, for $q>0$, we
have $y_2=e^{is}/\sqrt{q}$ for $s\in\R$, while for $q<0$, we have
$y_2=e^{is}j/\sqrt{-q}$ for $s\in\R$. In either case, $\cS_{0,q}\cap\HH2$
is identified with the set of $(iy_1,e^{is})\in \Imag\HQ\times S^1$ with
$y_1+\bar y_1<-p/q$. This is diffeomorphic to $S^1\times\R^3$.

It is now clear that when $\cS(p,q)\cap\HH2$ is non-empty, as in the previous
cases, it carries a complete SDE metric of negative scalar curvature. The
isometry group contains the $2$-torus
\begin{equation*}
(y_1,y_2)\mapsto (\sigma y_1\sigma^{-1}, \tau y_2\sigma^{-1})
\end{equation*}
with $(\tau,\sigma)\in U(1)\times U(1)$.
\end{proof}

\section{The Height Two Quotients}

To complete our analysis of the quotients of $\HH2$ we consider the height two
case $\tilde T_2(\lambda,ip)$. As in the height one case, by scaling and
conjugation, we can suppose $\Delta_2(p)=\tilde T_2(1,ip)$ with $p\in\{0,1\}$,
so there are only two distinct quotients up to scale, but we shall carry out
our computations for arbitrary $p$. We then have
\begin{equation}
\varphi_t^{p}(\bfv)=
\begin{pmatrix}
e^{ipt}&0&0\\ 
-t^2/2&e^{ipt}&it\\
it& 0 & e^{ip t}
\end{pmatrix}
\begin{pmatrix}v_0\\ v_1\\ v_2\end{pmatrix}
=\begin{pmatrix}e^{ipt}v_0\\ e^{ipt}v_1+it v_2-t^2v_0/2 \\ 
e^{ipt}v_2+it v_0\end{pmatrix},
\end{equation}
which reduces, in inhomogeneous coordinates $\bfy=(y_1,y_2)$, to
\begin{equation}\label{eq:h2-act}
\varphi_t^{p}\begin{pmatrix}y_1\\ y_2\end{pmatrix}=
\begin{pmatrix}
e^{ip t} (y_1+it y_2-t^2/2) e^{-ip t}\\
e^{ip t} (y_2+it) e^{-ipt}
\end{pmatrix}.
\end{equation}
The moment map for this action is given in homogeneous or inhomogeneous
coordinates by
\begin{align*}
\mu_{p}(\bfv) &= {\bar v_0} i v_2 +{\bar v_2}i v_0
+p({\bar v_0}iv_1+{\bar v_1}iv_0) + p {\bar v_2}i v_2,\\
f_{p}(\bfy)&= i y_2 +{\bar y_2}i +p(iy_1+{\bar y_1}i) + p {\bar y_2}i y_2.
\end{align*}

\begin{theorem} Let $\Delta= \Delta_2(p)=\tilde T_2(1,ip)$ and consider the
one parameter group $\varphi^{p}_t=e^{\Delta t}$ acting on the quaternionic
hyperbolic space $\HH2$. Then the quaternion K\"ahler reduction
$M(p)=\mu^{-1}_{p}(0)/\varphi^{p}$ is diffeomorphic to $\R^4$ for all $p$, and
carries a complete self-dual Einstein metric of negative scalar curvature whose
isometry group contains $S^1\times\R$. Furthermore, $M(0)$ is quaternionic
hyperbolic space.
\end{theorem}

\begin{proof} Consider the following set
\begin{equation}
\cS_p=\{\bfy\ \ |\ \ f_p(y_1,y_2)=
0\text{ and } iy_2-\bar y_2 i = 2\Real(iy_2)=0 \}.
\end{equation}
It is clear that this is a global slice for the action of $e^{\Delta t}$ on
the zero set of the momentum map. For $p=0$, we obtain $y_2=0$, and hence
$\cS_p\cap\HH2$ is diffeomorphic to $\{y_1\in\HQ:\Real y_1<0\}$, so let us
suppose that $p\neq 0$. We write $y_2=s_2+jw_2$ with $s_2\in\R$ and
$w_2\in\C$. The momentum constraint determines the imaginary part of $iy_1$ in
terms of $y_2$. In particular, it implies that
\begin{equation*}
2s_2/p + y_1+\bar y_1 + s_2^2-|w_2|^2 =0.
\end{equation*}
We find that $y_2$ is constrained to lie in the paraboloid $s_2/p >
|w_2|^2$. $M(p)$ is diffeomorphic to the product of this paraboloid with the
real line, which is diffeomorphic to $\R^4$, and as before has a complete
SDE metric of negative scalar curvature.

The isometry group of the quotient metric contains the group generated by
$\tilde T_1(\lambda,ip,ip)$, which is isomorphic to $S^1\times \R$.
The last statement follows from a direct computation.
\end{proof}

\section{The Bergman Metric on the 4-Ball}\label{berg}

In this section we turn our attention to the quotients of
$\HH{1,1}=\Proj_\HQ(\HQ^{1,2}_+)$.  One could consider all the cases studied
in the previous four sections. Locally we will get families of metrics of both
positive and negative scalar curvature. However, because $\HH{1,1}$ is not
Riemannian, singularities can arise when the vector field generating the
$\varphi_\Delta(t)=e^{\Delta t}$ action is null somewhere on the zero-set of
the momentum map. For this reason, we shall restrict our attention to the
special case $\Delta=\Delta_0(\bfp)$ (cf. \cite{MR88f:53088}). Furthermore, it
will be convenient to switch signature and take quotients of
$\Proj_\HQ(\HQ^{2,1}_-)$: this means we don't have to reverse the sign of the
quotient metric to get a positive definite metric of negative scalar
curvature.

We begin by placing the case $\bfp=(1,1,1)$ in a more general context. Recall
the following construction of the Wolf space $X(2,k)=U(2,k)/U(2)\times
U(k)$. We start with the space $\Proj_\HQ(\HQ^{2,k}_-)$ and the diagonal
circle action on $\Proj_\HQ(\HQ^{2,k}_-)$, described in quaternionic
coordinates $\bfu=(u_0,u_1,u_2,\ldots,u_{k+1})$ as
\begin{equation}\varphi_t(\bfu)=e^{2\pi it}\bfu,
\end{equation}
where $t\in[0,1/2)$. The moment map for this action reads
\begin{equation}\mu(\bfu)=-\bar u_0iu_0-\bar u_1iu_1+
\sum_{\alpha=2}^{k+1}\bar u_\alpha iu_\alpha.
\end{equation}
By introducing the complex coordinates $u_\alpha=z_\alpha+j\bar w_\alpha$
and the matrices
\begin{equation}
\Z=\begin{pmatrix}\Z_0\\ \Z_1\end{pmatrix},\ \ \ \ 
\Z_0=\begin{pmatrix}z_0&w_0\\ z_1&w_1\end{pmatrix},\ \ \ \ \
\Z_1=\begin{pmatrix}
z_2&w_2\\ \vdots&\vdots\\ z_{k+1}&w_{k+1}\end{pmatrix},
\end{equation}
we can describe the set $\mu^{-1}(0)\cap \cH_{2,k}(-1)$ by a matrix equation
\begin{equation}-\Z_0^\dagger\Z_0+\Z_1^\dagger\Z_1=-\mathbb I_{2\times 2}.
\end{equation}
Now, one observes that the $U(1)\cdot \Symp(1)\simeq U(2)$ which takes us from
$\mu^{-1}(0)\cap \cH_{2,k}(-1)$ to the quotient is nothing but $U(2)$ matrix
multiplication of $\Z$ from the right. This action is free and the quotient is
simple a bounded domain in $\C^{2k}$. As homogeneous (symmetric) spaces
\begin{equation}\mu^{-1}(0)\cap \cH_{2,k}(-1)\simeq U(2,k)/U(k),
\end{equation}
and 
\begin{equation}M= \frac{\mu^{-1}(0)\cap \cH_{2,k}(-1)}{U(2)}=
\frac{U(2,k)}{U(2)\times U(k)}.\end{equation}

In particular, when $k=1$ we get the complex hyperbolic (or Bergman) metric on
the until ball in $\C^2$.

Below, we will show that this construction is rigid in a sense that an
introduction of weights automatically leads to orbifold singularities.  As we
are interested in 4-dimensional quotients we will set $k=1$.  In the previous
sections we have seen that all of the complete $U(2)$-symmetric SDE metrics of
negative scalar curvature can be obtained as quaternion K\"ahler quotients of
the ball $\Proj_\HQ(\HQ_-^{1,2})$.  The only exception is the complex
hyperbolic Bergman metric. The above calculation now shows that this metric
can be constructed as a quotient of the pseudo-Riemannian quaternion K\"ahler
manifold $\Proj_\HQ(\HQ_-^{2,1})$. More generally, take $\Delta=\Delta(\bfp)$
and examine the following circle action
\begin{equation}\varphi^\bfp_t(u_0,u_1,u_2)=
(e^{2\pi ip_0t}u_0,\ e^{2\pi ip_1t}u_1,
\ e^{2\pi ip_2t}u_2)\end{equation}
where $\bfp=(p_0,p_1,p_2)\in\Z^3$,
${\rm gcd}(p_0,p_1,p_2)=1$, $t\in[0,1)$ when all the weights are odd,
and $t\in[0,1/2)$ otherwise. Now,
the moment map $\mu_\bfp:\Proj_\HQ(\HQ_-^{2,1})\rightarrow\cV$ 
is given as
\begin{equation}
\mu_\bfp(\bfu)=-p_0\bar u_0iu_0-p_1\bar u_1iu_1+p_2\bar u_2iu_2.
\end{equation}

\begin{theorem} Let $\bfp\in(\Z^+)^3$ and let $M(\bfp)$ be the
quaternion K\"ahler quotient of $\Proj_\HQ(\HQ^{2,1}_-)$ by the above circle
action. Then $M(\bfp)$ has orbifold singularities unless $\bfp=(1,1,1)$ in
which case $M(1,1,1)\simeq U(2,1)/U(2)\times U(1)$ is the symmetric complex
hyperbolic metric on the unit ball in $\C^2$.
\end{theorem}

\begin{proof}
Unlike in the case of $\Proj_\HQ(\HQ_-^{2,1})$ we no longer have the advantage
of global coordinates. We need to consider two cases
\begin{equation}\Proj_\HQ(\HQ_-^{2,1})=\cU_0\cup \cU_1,\end{equation}
where $\cU_i$ are defined as a submanifold on which $u_i\not=0$.
We first consider $\cU_0$, where we can switch to inhomogeneous local chart
$(x^0_1,x^0_2)=(u_1u_0^{-1}, u_2u_0^{-1})$. On $\cU_0$ we have
\begin{equation}\label{7.10}
-|x^0_1|^2+|x^0_2|^2<1.
\end{equation}
As before, the action and the zero level of the moment map become
\begin{equation}\varphi^\bfp_t(x^0_1,x^0_2)=
(e^{2\pi ip_1t}x^0_1 e^{-2\pi ip_0t},
e^{2\pi ip_2t}x^0_2 e^{-2\pi ip_0t}),\end{equation}
\begin{equation}\label{7.12}
0=-ip_0-p_1\bar x^0_1ix^0_1+p_2\bar x^0_2ix^0_2.
\end{equation}
We then write
\begin{equation}
\bfx^0=\bfz^0+j\bar \bfw^0,
\end{equation}
where $(\bfz^0,\bfw^0)\in \cU_0$ 
and observe that on $\cU_0$
\begin{equation*}
\varphi^\bfp_t\begin{pmatrix}z^0_1&w^0_1\\ z^0_2&w^0_2\end{pmatrix}=
\begin{pmatrix}e^{2\pi i(p_1-p_0)t}z^0_1&e^{2\pi i(p_1+p_0)t}w^0_1\\
e^{2\pi i(p_2-p_0)t}z^0_2&e^{2\pi i(p_2+p_0)t}w^0_2\end{pmatrix}
\end{equation*}
while the moment map equations~\eqref{7.12} become 
\begin{equation}
-p_1(|z^0_1|^2-|w^0_1|^2)+p_2(|z^0_2|^2-|w^0_2|^2)=p_0,\ \ \ 
-p_1\bar w^0_1z^0_1+ p_2\bar w^0_2z^0_2=
0.\end{equation}
In this case we no longer have any analogue of Proposition~\ref{prop3}
as $\mu_\bfp^{-1}(0)$ always intersects the open set defined by~\eqref{7.10}.

Without loss of generality we will further assume that all weights are
non-negative. Furthermore, neither $p_0$ nor $p_1$ can equal 0 if we want the
quotient to be non-singular. If, say, $p_0=0$ then take $(u_0,0,0)\in
\Proj_\HQ(\HQ^{2,1}_-)$.  This point is also on the level set of the moment
map and it is fixed by every element of $S^1(\bfp)$. The third weight $p_2$
can be zero.  On $\cU_0\cap \mu_\bfp^{-1}(0)$ we can choose
$z^0_1=z^0_2=w^0_2=0$ and $|w^0_1|^2=p_0/p_1$ which is a circle of points
fixed by $\Z_{p_0+p_1}$. Hence, in order to get smooth quotient we must assume
all $p_0=p_1=1$ and $p_2$ is odd.  But then, taking $z^0_1=w^0_1=w^0_2=0$ and
$|z^0_2|^2=p_0/p_2=1/p_2$ one gets a circle of points where the isotropy group
equals $\Z_{(p_2+p_0)/2}$.  This forces $p_0=p_1=p_2=1$.  From our
previous example we know that $M(1,1,1)$ is the complex hyperbolic Bergman
metric on $\C^2$.
\end{proof}

\begin{remark} Let us observe that all quaternion K\"ahler reductions of
the symmetric space $X(2,2)\simeq U(2,2)/U(2)\times U(2)$ can now be obtained
using our construction in a very simple manner. As $X(2,2)$ is by itself
reduction of $\Proj_\HQ(\HQ^{2,2}_-)$ by the circle action corresponding to
the generator $T_1=ip\bbi_4$ we can consider all possible quotients of
$\Proj_\HQ(\HQ^{2,2}_-)$ by $2$-dimensional Lie algebras $\g=\{T_1, T_2\}$,
where $T_2\in\symp(2,2)$. Since $T_1$ is fixed to be a multiple of the
identity these are classified by the adjoint orbit $[T_2]$ in
$\symp(2,2)$. Hence, one could begin by enumerating all such classes. Here,
there are many more cases. To begin with $\symp(2,2)$ has 3 different Cartan
subalgebras. In addition, we have elements of height 0,1,2,3. In fact, there
are six distinct families of `purely' nilpotent classes \cite{MR55:5761}.
All of these quotients can be examined and they should lead to many new
metrics.
\end{remark}

\begin{example} The simplest example when one gets non-trivial negative
SDE Hermitian metric is deformation of the Bergman metric on the 4-ball. We
choose the second generator as
\begin{equation}
T_2(p,q,r)=
\begin{pmatrix}
ip&0&0&0\\ 
0&iq&1&0\\ 
0&1&iq&0\\ 
0&0&0&ir 
\end{pmatrix} \in \symp(2,2),
\end{equation}
One can easily see that $p=q=r=0$ gives the Bergman metric which should
correspond to a 4-parameter family of deformations of this metric. Detailed
analysis of this and other quotients will be carried out elsewhere.
\end{example}

\section{Quotients, hyperbolic eigenfunctions and Bochner-flat metrics}

The SDE metrics that we have constructed have in common that they possess (at
least) two commuting Killing vector fields, and therefore belong to the class
of metrics classified locally by Calderbank and Pedersen~\cite{MR1950174}.
Furthermore, according to Apostolov--Gauduchon~\cite{math.DG/0003162},
quaternion K\"ahler quotients of $\HH2$ are not just SDE, but Hermitian, and
are therefore conformal to the self-dual (and therefore Bochner-flat) K\"ahler
metrics classified by Bryant~\cite{MR2002i:53096}. In this section we relate
our metrics to the hyperbolic eigenfunction Ansatz of Calderbank--Pedersen
(which gives the explicit local form of the metrics), and to the SDE Hermitian
metrics of Apostolov--Gauduchon and Bryant.

We recall that the work of Calderbank and Pedersen shows that an SDE metric of
nonzero scalar curvature with two commuting Killing vector fields is
determined explicitly (on the open set where the vector fields are linearly
independent) by an eigenfunction $F$ of the Laplacian on the hyperbolic plane
with eigenvalue $3/4$, so it suffices to find the eigenfunction $F$
corresponding to our quotients. According to~\cite{MR1950174}, the hyperbolic
eigenfunctions $F$ arising as quotients of $\HH2$ or $\HH{1,1}$ should be
either `$3$-pole' solutions, or limits in which one or more of the `centers'
of the $3$-pole coincide. We shall justify this claim here.

\subsection{Quotients and hyperbolic eigenfunctions}

We first consider SDE manifolds arising as semi-quaternion K\"ahler quotients
of $\HH{k-1,l}$ or $\HH{k,l-1}$ by $n-1$ dimensional Abelian subgroups $G$ of
$\Symp(k,l)$ (with $k+l=n+1$) in full generality. Following~\cite{MR1950174},
we study such a quotient $(M,g)$ using the Swann bundle $(\tilde M,\tilde g)$,
which is the principal $\CO(3)$ bundle over $(M,g)$ arising as the
corresponding semi-hyperk\"ahler quotient of $\HQ^{k,l}$. More precisely, we
take the semi-hyperk\"ahler quotient by $G$ of (a connected component of)
$\HQ^{k,l}_*=(\HQ^{k,l}\smallsetminus\HQ^{k,l}_0)/\{\pm1\}$, which is a
principal $\CO(3)$-bundle over $\HH{k-1,l}\cup\HH{k,l-1}$.  $\tilde M$ is thus
the quotient by $G$ of the zero-set of the momentum map of $G$ in
$\HQ^{k,l}_*$ and we have a commutative diagram
\begin{equation*}
\begin{matrix}
\HQ^{k,l}_* &&\longrightarrow&& \HH{k-1,l}\cup\HH{k,l-1} \\
 \Big\downarrow &&&&\Big\downarrow\\
\tilde M &&\longrightarrow&& M,
\end{matrix}
\end{equation*}
where the vertical arrows denote semi-hyperk\"ahler and quaternion K\"ahler
quotients, and the horizontal arrows are principal $\CO(3)$ bundles:
$\SO(3)$ acts by isometries, and $\R^+$ by homotheties, so that if $q$ is
an $\HQ^\times$ valued function on the double cover of $\tilde M$ coming
from a local trivialization, we have
\begin{equation*}
\tilde g = s |q|^2 g + |dq +q\omega|^2,
\end{equation*}
where $\omega$ is the principal $\SO(3)$ connection on $\tilde M$ and $s$ is a
positive multiple of the scalar curvature of $g$, so that $s g$ is a (possibly
negative definite) SDE metric of positive scalar curvature. We can then
arrange our conventions so that $|q|^2$ pulls back to the zero-set of the
momentum map in $\HQ^{k,l}$ to give the absolute value $|F_{k,l}(\bfu,\bfu)|$
of the quadratic form.

Any $(n-1)$-dimensional Abelian subgroup $G$ of $\Symp(k,l)$ lies in a maximal
Abelian subgroup $H$, which has dimension $n+1$. For generic $G$ this maximal
Abelian subgroup will be unique, but in general we must choose such a group
$H$ so that we have a quotient group $H/G$ acting on $\tilde M$ and $M$. The
Lie algebra of this quotient group may be identified with $\R^2$.

Now, according to~\cite{MR1950174}, there is also a commutative diagram
\begin{equation*}
\begin{matrix}
\tilde M &&\longrightarrow&& M\\
 \Big\downarrow &&&&\Big\downarrow\\
\Imag \HQ\tens_0 \R^2 &&\longrightarrow&& \cH^2,
\end{matrix}
\end{equation*}
where the vertical arrows are (possibly only locally defined) isometric
quotients by $H/G$, $\Imag \HQ\tens_0\R^2$ is the open subset of
indecomposable elements of the flat vector space $\Imag \HQ\tens
\R^2\cong\Imag\HQ\dsum\Imag\HQ$, and $\cH^2$ is the hyperboloid of positive
definite elements of determinant one in $S^2\R^2$, equipped with the metric
induced by the determinant on $S^2\R^2$ (which is the hyperbolic metric). The
lower horizontal map, like the upper map, is a principal $\CO(3)$-bundle, and
is given explicitly by the Grammian map
\begin{equation*}
\bfx=(x_1,x_2)\in \Imag\HQ \tens_0 \R^2\to \frac1{|x_1\wedge x_2|}
\begin{pmatrix}
|x_1|^2 & \ip{x_1,x_2}\\
\ip{x_1,x_2} & |x_2|^2
\end{pmatrix}.
\end{equation*}
Given a hyperbolic eigenfunction $F$ on (an open subset of) $\cH^2$, we can
lift $F$ to a homogeneity $1/2$ function $\tilde F$ on the corresponding union
of rays in the space of positive definite elements of $S^2\R^2$. Now we have
the following result, which was proven in the definite case (i.e., $k=0$ or
$l=0$) in~\cite{MR1950174}. The more general result also has a more direct
proof, and we correct a minor error in~\cite{MR1950174}.

\begin{theorem} Let $(M^4,g)$ be an SDE metric with two commuting Killing
vector fields obtained as a semi-quaternion K\"ahler quotient of $\HH{k-1,l}$
or $\HH{k,l-1}$ by an $n-1$ dimensional Abelian subgroup $G$ of $\Symp(k,l)$
\textup(where $k+l=n+1$\textup), and let $\tilde F$ be the homogeneity $1/2$
lift of the hyperbolic eigenfunction $F$ generating $g$, locally, with respect
to a $2$-dimensional Abelian quotient group acting by isometries. Then the
pullback of $\tilde F$ to the zero-set of the momentum map in $\HQ^{k,l}$ is
a nonzero constant multiple of the restriction of the quadratic form
$F_{k,l}(\bfu,\bfu)$.
\end{theorem}
\begin{proof}
Let $A$ be a positive definite element of $S^2\R^2$, and write $A=\sqrt{\det
A} A_1$ with $A_1\in\cH^2$ having determinant one. Then by definition $\tilde
F(A) = (\det A)^{1/4} F(A_1)$ and so $\tilde F = (\det A)^{1/4} F$, where $F$
now denotes the (homogeneity $0$) pullback to $S^2\R^2$. Now it was shown
in~\cite{MR1950174} that the pullback of the function $A\mapsto\det A$ to the
Swann bundle $\tilde M$ is $|q|^8/|F|^4$ (although the result is incorrectly
stated there). It follows that $\tilde F$ pulls back to the Swann bundle to
give $|q|^2 F/|F|$, which pulls back to the momentum zero set in $\HQ^{k,l}$
to give a nonzero constant multiple of the absolute value of the quadratic
form times a (possibly nonconstant) sign.  However, $\tilde F$ is smooth, even
through its zero-set, so the result follows.
\end{proof}

Note that the pullback of $\tilde F$ is independent of the choice of quotient
torus (in the case that such a choice exists).

We are going to use this result to calculate the hyperbolic eigenfunction
corresponding to the metrics we have studied in detail here. In order to do
this we just need to write the quadratic form $F_{k,l}(\bfu,\bfu)$ in momentum
coordinates and restrict it to the zero-set of the momentum map, as we now
explain.

Having chosen (if there is a choice) the maximal Abelian subgroup $H$ of
$\Symp(k,l)$ containing $G$ (and a basis for the Lie algebra of $H$ so that we
can identify it with $\R^{n+1}$), we have momentum coordinates $y_0,\ldots
y_n\in \Imag\HQ$ which are independent on the open subset $U$ of $\HQ^{k,l}$
where the $H$ action is free. Since $F_{k,l}(\bfu,\bfu)$ is $H$-invariant it
will be a function of $\bfy=(y_0,\ldots y_n)$ on $U$, and our first task is to
compute this function. Then, secondly, we must restrict to the zero-set of the
momentum map of the $G$ action.  For this second step, following
Bielawski--Dancer~\cite{BiDa:gtt}, we introduce an explicit parameterization
of the momentum zero-set of $G$ in terms of the momentum coordinates of the
quotient torus. To do this, we write the Lie algebra $\g$ of $G$ as the kernel
of a $2\times (n+1)$ matrix $S\colon \R^{n+1} \to\R^2$. Then the transpose
matrix $S^t\colon \R^{2*}\to \R^{(n+1)*}$ parameterizes the kernel of the
projection $\R^{(n+1)*}\to\g^*$. Since the momentum map of $H$ is injective on
$U$, the momentum zero-set of $G$ in $U$ is the subset where the momentum map
of $H$ takes values in the image of $S^t$, so we can parameterize it by
writing $\bfw=S^t\bfx$, with $\bfx=(x_1,x_2)$.

The hyperbolic eigenfunction is now obtained by substituting this into the
quadratic form $F_{k,l}$, writing the result in terms of the $\SO(3)$
invariants $\ip{x_i,x_j}$ and restricting to the hyperboloid $\det
\ip{x_i,x_j} = 1$, where we can write
\begin{equation*}
\begin{pmatrix}
|x_1|^2 & \ip{x_1,x_2}\\
\ip{x_1,x_2} & |x_2|^2
\end{pmatrix}=
\begin{pmatrix}
1/\rho &\eta/\rho\\
\eta/\rho &(\rho^2+\eta^2)/\rho
\end{pmatrix}
\end{equation*}
for half-space coordinates $(\rho,\eta)$ on $\cH^2$.  We now carry out this
procedure for the examples we have studied.

\subsection{Subgroups of a maximal torus and the generalized Pedersen--LeBrun
metrics}

Let $H\cong (S^1)^{n+1}$ be the standard maximal torus in $\Symp(k,l)$ acting
diagonally on $\HQ^{k,l}$ with respect to the coordinates $(u_0,\ldots u_n)$,
i.e., the $j$th circle acts by scalar multiplication by $e^{it}$ on the $j$th
coordinate $u_j$, and has momentum map $y_j = \bar u_j i u_j$. We
therefore have
\begin{equation*}
F_{k,l}(\bfu,\bfu) = -\sum_{j=0}^{k-1} |y_j|+\sum_{j=k}^{k+l}|y_j|.
\end{equation*}
On the zero-set of the momentum map of $G$ we then get
\begin{equation*}
F_{k,l}(\bfu,\bfu) = -\sum_{j=0}^{k-1} |a_j x_2 - b_j x_1|
+\sum_{j=k}^{k+l} |a_j x_2 - b_j x_1|,
\end{equation*}
where the matrix $S_{ij}$ defining $\g$ has columns $(-b_j,a_j)$.
We now observe that
\begin{align}\notag
|a x_2 - b x_1|
&= \sqrt{ a^2 |x_2|^2 - 2 a b \ip{x_1,x_2} + b^2 |x_1|^2 }\\
&= \frac{\sqrt{ a^2 (\rho^2+\eta^2) - 2 a b \eta + b^2 }}{\sqrt{\rho}}
= \frac{\sqrt{ a^2 \rho^2 + (a \eta - b)^2 }}{\sqrt{\rho}}
\label{monopole}
\end{align}
and thus the corresponding hyperbolic eigenfunction is
\begin{equation*}
F(\rho,\eta) = -\sum_{j=0}^{k-1} \frac{\sqrt{ a_j^2\rho^2 +
(a_j\eta-b_j)^2 }}{\sqrt{\rho}}+\sum_{j=k}^{k+l} \frac{\sqrt{ a_j^2\rho^2 +
(a_j\eta-b_j)^2 }}{\sqrt{\rho}}
\end{equation*}
in accordance with the discussion in \cite{MR1950174}---see
also~\cite{CaSi:emcs}.

These hyperbolic eigenfunctions may be interpreted as `multipole' solutions,
in the sense that they are a linear combination of solutions of the
form~\eqref{monopole} which we regard as the eigenfunction generated by a
monopole source at the point $\eta=b/a$ on the boundary $\rho=0$ of the
hyperbolic plane (which is a circle $\R\cup\{\infty\}$).

In the case studied in this paper, $n=3$, and the vectors $(a_0,a_1,a_2)$ and
$(b_0,b_1,b_2)$ are any two linearly independent solutions to the equation
$a_0p_0+a_1p_1+a_2p_2=0$, where $p_0,p_1,p_2$ are the weights of the torus
action.

Note that $\SL(2,\R)$ acts on the vectors $(a_j,b_j)$ to give equivalent
solutions so that the points $b_j/a_j$ can be fixed (for instance at $1,-1$
and $\infty$, as in \cite{MR1950174}---we remark that the points are distinct
provided the weights $p_0,p_1,p_2$ are nonzero).

\subsection{The generalized Pedersen metrics}

For the generalized Pedersen metrics, the family of generators that we are
using span a Cartan subalgebra of $\symp(1,2)$ which is not the Lie algebra of
a maximal torus. However, this can be understood as an analytic continuation
of the generalized Pedersen--LeBrun metrics (replace $\lambda$ by $it$, where
$i$ acting on the left is a complex scalar commuting with the right
quaternionic action, and diagonalize). As discussed in~\cite{MR1950174} this
implies that the hyperbolic eigenfunction can be assumed to take the form
\begin{equation*}
F(\rho,\eta)= \frac a{\sqrt\rho}
+\frac{b+ic}2\frac{\sqrt{\rho^2+(\eta+i)^2}}{\sqrt\rho}
+\frac{b-ic}2\frac{\sqrt{\rho^2+(\eta-i)^2}}{\sqrt\rho}.
\end{equation*}
This is still a $3$-pole solution, but two of the sources are complex
conjugate rather than real.

\subsection{The height one quotients}

In the remaining cases, it is more convenient to begin with the coordinates
$\bfv=(v_0,v_1,v_2)$ that we introduced already before, so that
\begin{equation*}
F_{1,2}(\bfu,\bfu) = \bar v_0 v_1+ \bar v_1 v_0 + \bar v_2 v_2.
\end{equation*}

In the case of the height one quotients, the momentum coordinates (in terms of
the $v_0,v_1,v_2$ coordinates) that we shall use are
\begin{equation*}
y_0 = \bar v_1 i v_0 + \bar v_0 i v_1, \qquad
y_1 = -\bar v_0 i v_0,\qquad
y_2 = \bar v_2 i v_2
\end{equation*}
and we compute
\begin{equation*}
F_{1,2}(\bfu,\bfu) = \bar v_0 v_1+ \bar v_1 v_0 + \bar v_2
v_2 = \frac{\ip{y_0,y_1}}{|y_1|} + |y_2|.
\end{equation*}
After substituting for $x_1,x_2$, the second term is treated as before, so it
suffices to compute
\begin{equation*}
\frac{\ip{a_0 x_2-b_0 x_1, a_1 x_2-b_1 x_1 }}
{|a_1 x_2-b_1 x_1|}
= \frac{ a_0 a_1\rho^2 +(a_0\eta-b_0)(a_1\eta-b_1)}
{\sqrt{\rho}\sqrt{a_1^2\rho^2+(a_1\eta-b_1)^2} }.
\end{equation*}

Under the action of $\SL(2,\R)$ this is equivalent to
\begin{equation*}
F(\rho,\eta)= \frac{\eta}{\sqrt\rho\sqrt{\rho^2+\eta^2}}
=\frac{\partial}{\partial\eta}\frac{\sqrt{\rho^2+\eta^2}}{\sqrt\rho}
\end{equation*}
which may be interpreted as an `infinitesimal dipole', i.e., a limit of
oppositely charged monopoles at $\eta=\pm\eps$ as $\eps\to0$.

Thus the hyperbolic eigenfunctions corresponding to height one quotients are
combinations of a monopole and an infinitesimal dipole.

\subsection{The height two quotients}

In the height two case, the maximal Abelian subalgebra containing
$T_2$ is unique, being spanned by $i\bbi_3$, $T_2$ and $T_1=i T_2^2$.
The momentum coordinates of these generators are
\begin{equation*}
y_0 = \bar v_1 i v_0 + \bar v_0iv_1 + \bar v_2iv_2,\qquad
y_1 = \bar v_2 i v_0 + \bar v_0 i v_2, \qquad
y_2 = -\bar v_0 i v_0.
\end{equation*}
Writing the quadratic form in these coordinates is straightforward once
one has computed all the inner products between them. The result is
\begin{equation*}
F_{1,2}(\bfu,\bfu) = \bar v_0 v_1+ \bar v_1 v_0 + \bar v_2
v_2 = \frac{|y_1|^2 |y_2|^2 - \ip{y_1,y_2}^2 +
2\ip{y_0,y_2}|y_2|^2} {2|y_2|^3}.
\end{equation*}
The family of quotients we consider is the span of $i\bbi_3$ and $T_2$, so
we can take $y_0=a_0 x_1$, $y_1=a_1 x_1$ and $y_2=x_2$
as our parameterization in quotient coordinates to yield
\begin{equation*}
\frac{a_1^2 (|x_1|^2 |x_2|^2 - \ip{x_1,x_2}^2) +
2a_0 \ip{x_1,x_2}|x_2|^2} {2|x_2|^3}=
a_0 \frac {\eta}{\sqrt{\rho}\sqrt{\rho^2+\eta^2}}+\frac{a_1^2}{2}
\frac{\rho^{3/2}}{(\rho^2+\eta^2)^{3/2}}.
\end{equation*}
We recognise the first term as an infinitesimal dipole. Differentiating again
with respect to $\eta$, we see that the second term may be regarded as an
infinitesimal tripole. As the two terms have different homogeneities in
$(\rho,\eta)$, by scaling the coordinates and the eigenfunction, we have
just three distinct quotients:
\begin{itemize}
\item $a_0=0$, the pure tripole, corresponds the quotient by $i\bbi_3$, which
is complex hyperbolic space (under the non-semisimple $\R^2$ action induced by
$T_1$ and $T_2$);
\item $a_1=0$, the pure dipole, corresponds to the quotient by $T_2$, which
is real hyperbolic space (under the non-semisimple $S^1\times\R$ action
induced by $T_1$ and $i\bbi_3$);
\item the nontrivial case with $a_0,a_1$ both nonzero.
\end{itemize}

\subsection{Infinitesimal multipoles from the quotient point of view}

We have seen, as conjectured in \cite{MR1950174}, that the nilpotent cases
(height one and height two quotients) can be regarded as limiting cases in
which two or more monopoles come together to form an infinitesimal
multipole. This can be seen from the group theory of the quotient construction
by realizing a non-semisimple element as a limit of semisimple ones.

For example, consider the following generator in $\symp(1,2)$:
\begin{equation}
\mathbb T_{\bfp,\lambda}=\begin{pmatrix}
ip_0&\lambda&0\\ \lambda&ip_1&0\\ 0&0&ip_2
\end{pmatrix} \in \symp(1,2).\end{equation}
Generically this generator is of height 0 but a special choice of the
parameters $(\bfp,\lambda)$ raises the height to 1. To see it let us consider
the following one parameter group actions on the ball
\begin{equation*}
\varphi_t^{\bfp,\lambda}(\bfu) =\exp(\mathbb
T_{\bfp,\lambda}t)\cdot\bfu\equiv \mathbb A_{\bfp,\lambda}(t)\cdot\bfu ,
\end{equation*}
where now all $(\lambda,p_0,p_1,p_2)$ are real parameters and $\mathbb
A_{\bfp,\lambda}(t)\in U(1,2)\subset Sp(1,2)$ and we assume $\lambda\neq 0$.
We set
\begin{equation}\alpha=\frac{p_0-p_1}{2},\ \ \ \beta=\frac{p_0+p_1}{2},\ \ \ \ 
\gamma=\sqrt{|\alpha^2- \lambda^2|}.
\end{equation}
We can compute the matrix $\mathbb A_{\bfp,\lambda}(t)$ explicitly.  Depending
on the sign of $\alpha^2-\lambda^2$ we get three distinct cases.  If we denote
the corresponding $U(1,2)$ matrices by $\mathbb A_{\bfp,\lambda}^+(t), \mathbb
A_{\bfp,\lambda}^0(t),$ and $\mathbb A_{\bfp,\lambda}^-(t)$, we obtain
\begin{align*}
\mathbb A_{\bfp,\lambda}^+(t)&=\begin{pmatrix}
e^{i\beta t}(\cosh\gamma t+\frac{i\alpha}{\gamma}\sinh\gamma t)&
\frac{\lambda}{\gamma}e^{i\beta t}\sinh\gamma t&0\\
\frac{\lambda}{\gamma}e^{i\beta t}\sinh\gamma t&
e^{i\beta t}(\cosh\gamma t-\frac{i\alpha}{\gamma}\sinh\gamma t)&0\\
0&0&e^{ip_2t}\end{pmatrix},\\
\mathbb A_{\bfp,\lambda}^0(t)&=\begin{pmatrix} 
e^{i\beta t}(1+i\alpha t)& 
e^{i\beta t}\lambda t&0\\ 
e^{i\beta t}\lambda t& 
e^{i\beta t}(1-i\alpha t)&0\\ 
0&0&e^{ip_2t}\end{pmatrix},\\
\mathbb A_{\bfp,\lambda}^-(t)&=\begin{pmatrix}
e^{i\beta t}(\cos\gamma t+\frac{i\alpha}{\gamma}\sin\gamma t)&
\frac{\lambda}{\gamma}e^{i\beta t}\sin\gamma t&0\\
\frac{\lambda}{\gamma}e^{i\beta t}\sin\gamma t&
e^{i\beta t}(\cos\gamma t-\frac{i\alpha}{\gamma}\sin\gamma t)&0\\
0&0&e^{ip_2t}\end{pmatrix}.
\end{align*} 
Note that $\lim_{\gamma\to0}\mathbb A_{\bfp,\lambda}^+(t)=
\lim_{\gamma\to0}\mathbb A_{\bfp,\lambda}^-(t)=\mathbb A_{\bfp,\lambda}^0(t)$.
Also, $\mathbb A_{\bfp,\lambda}^-(t)$ is actually a circle provided the triple
$(\gamma,\beta,p_2)$ is commensurate (all ratios are in $\Q$).

Note that the above calculation has to do with writing
\begin{equation}
T_{\bfp,\lambda}=L+N=
\begin{pmatrix}
i\beta&0&0\\ 0&i\beta&0\\ 0&0&ip_2
\end{pmatrix} 
+
\begin{pmatrix}
i\alpha&\lambda&0\\ \lambda&-i\alpha&0\\ 0&0&0
\end{pmatrix},
\end{equation}
where $[L,N]=0$ and $L=T_0(i\beta,i\beta,ip_2)$. Now, $N^2=0$ when
$\lambda^2=\alpha^2$. This shows that when $\lambda^2=\alpha^2$,
$T_{\bfp,\lambda}$ must be conjugated to some $T_1(\mu,p,q)$ of
Definition~\ref{def2}. When $\lambda^2\not=\alpha^2$ the generator
$T_{\bfp,\lambda}$ has height 0 and, depending on the sign of
$\lambda^2-\alpha^2$, is conjugated either to some $T_0(\mu,ip,iq)$ or
$T_0(iq_0,iq_1,iq_2)$. In either case, one can think of height 1 metrics as
certain limits of height 0 metrics.

\subsection{Quotients and Bochner-flat K\"ahler metrics}

We finally discuss the relationship between quaternion K\"ahler quotients of
$\HH2$ or $\HH{1,1}$ and Bochner-flat (i.e., self-dual) K\"ahler surfaces. On
a self-dual K\"ahler surface $(M,h,J)$ the conformal metric $g=s_h^{-2} h$,
defined wherever the scalar curvature $s_h$ of $h$ is nonzero is an SDE
metric. Conversely, $h$ can be recovered from $g$ using the fact that the Weyl
tensor $W=W^+$ of a self-dual K\"ahler surface is a constant multiple of $s_h
\omega\tens_0\omega,$ where $\omega$ is the K\"ahler form and the subscript
zero denotes the tracefree part in $S^2_0(\Lambda^2T_+\dual M)$: thus, up to a
constant $s_h=|W|_h =|W|_g^{1/3}$ and $h=|W|^{2/3}g$. This sets up a one to
one correspondence, at least locally, between self-dual K\"ahler metrics and
SDE Hermitian metrics~\cite{MR84h:53060,math.DG/0003162} which are not
conformally flat.  ($h$ and $g$ are equal up to homothety iff they are locally
symmetric.)

Bochner-flat K\"ahler manifolds have been completely classified, locally and
globally, by Bryant~\cite{MR2002i:53096}. The local classification is quite
easy to understand: over a Bochner-flat K\"ahler $2n$-manifold $M$, the
(locally defined) rank $1$ bundle with connection, whose curvature is the
K\"ahler form of $M$, has a flat CR structure (given by the horizontal lift of
the K\"ahler structure on $M$) and is therefore locally CR isomorphic to
$S^{2n+1}$. This realises the K\"ahler metric on $M$ as local quotient of
$S^{2n+1}$ by a one parameter subgroup of $\PSU(1,n+1)$, the group of CR
automorphisms of $S^{2n+1}$ (which is naturally realised as the quadric of
totally null complex lines in the projective space of $\C^{1,n+1}$).  It then
follows that Bochner-flat K\"ahler metrics are classified by adjoint orbits
in $\su(1,n+1)$.

Specialising to $n=2$, self-dual K\"ahler surfaces are classified, as local
quotients of $S^5$, by adjoint orbits in $\su(1,3)$, and it is natural to
conjecture that the corresponding SDE Hermitian metrics are obtained as
(perhaps only local) quaternion K\"ahler quotients of $\HP2$, $\HH2$ and
$\HH{1,1}$, classified by adjoint orbits in $\symp(3)$ and $\symp(1,2)$.  This
is essentially correct, as the work of
Apostolov--Gauduchon~\cite{math.DG/0003162} shows.

\begin{proposition}
Let $(M,g)$ be a self-dual Einstein manifold given as a (semi-)quaternion
K\"ahler quotient of $\HP2$, $\HH2$ or $\HH{1,1}$ by a (possibly local) $S^1$
or $\R$ action. Then $(M,g)$ admits a compatible Hermitian structure and there
is an invariant Sasakian structure on the momentum zero-set of the action,
whose underlying CR structure is flat.
\end{proposition}
\begin{proof}[Sketch proof] We outline the arguments, refering the reader to
Apostolov--Gauduchon~\cite{math.DG/0003162} for more details.  Let be $K$ a
quaternionic Killing vector field on a (semi-)quaternion K\"ahler manifold $Q$
of nonzero scalar curvature; this means that
\begin{equation*}
\nabla K\in\Cinf(Q,\cV_Q\oplus\symp(TQ))\subset \Cinf(Q,\so(TQ)),
\end{equation*}
where $\cV_Q$ is the bundle of $\symp(1)$'s in $\so(TQ)$ defining the
quaternionic structure and $\symp(TQ)\subset\so(TQ)$ is the bundle of
$\symp(n)$'s in $\so(TQ)$ consisting the skew endomorphisms which commute with
$\cV_Q$. Since the scalar curvature is nonzero, the momentum map of $K$ is
defined to be the $\cV_Q$ component of $\nabla K$. It follows that on the
zero-set $\cS$ of the momentum map, $\nabla K$ is a section of $\symp(TQ)$. It
is also $K$ invariant, so its horizontal part descends to the (perhaps only
locally defined) quotient $M=\cS/K$, which is the quaternion K\"ahler quotient
of $Q$ by $K$, to give a section $\Psi$ of $\symp(TM)$. If $Q$ is an
$8$-manifold, then $M$ is a $4$-manifold and $\so(TM)=\cV_M\oplus\symp(TM)$
and with our conventions $\cV_M=\so_-(TM)$ and $\symp(TM)\cong\so_+(TM)$, the
bundles of (anti-)self-dual endomorphisms associated to $\Lambda^2_-T\dual M$
and $\Lambda^2_+T\dual M$ using the metric. It follows that wherever $\Psi$ is
nonzero $\sqrt2 \Psi/|\Psi|$ is a almost complex structure which is self-dual
(i.e., orthogonal and commuting with the quaternionic structure), so that $M$
is an almost Hermitian manifold.  Apostolov and Gauduchon show that this
complex structure is integrable if $Q$ is $\HP2$ or $\HH2$ and their argument
applies unchanged to $\HH{1,1}$ (it is a straightforward consequence of the
fact that these spaces are flat as quaternionic manifolds). Thus $M$ is SDE
Hermitian, as claimed, and one can check that the conformal K\"ahler metric
$h$ is $|K|^{-2} g$.

Now the curvature of the rank $1$ bundle $\cS\to M$ (i.e., the horizontal part
of the $2$-form associated to $|K|^{-2}\nabla K$) is then the K\"ahler form of
$M$, so that the K\"ahler structure on $M$ lifts to the horizontal
distribution to give a $K$-invariant Sasakian structure on $\cS$. This is the
canonical Sasakian structure associated to $(M,h)$, and the underlying CR
structure is flat because $(M,h)$ is self-dual.
\end{proof}

A flat CR manifold is locally isomorphic to $S^5$ with its standard flat CR
structure (as the projective light cone in $\C^{1,3}$). Since $K$ generates an
action by CR automorphisms, such a local isomorphism determines an element of
$\su(1,3)$, the Lie algebra of CR automorphisms of $S^5$. However, the local
isomorphism is only determined up to conjugation by $\PSU(1,3)$, so we do not
obtain a Lie algebra homomorphism from $\symp(3)$ or $\symp(1,2)$ to
$\su(1,3)$---these Lie algebras are certainly not isomorphic.

Nevertheless, the classifications of self-dual K\"ahler manifolds (in terms of
adjoint orbits in $\su(1,3)$) and quotients of $\HP2$, $\HH2$ and $\HH{1,1}$
(in terms of adjoint orbits in $\symp(3)$ and $\symp(1,2)$) do essentially
coincide. This is slightly subtle, as in both quotient constructions the
manifold (or orbifold) corresponding to a conjugacy class may not be
connected: for the K\"ahler metric, these components correspond to Bryant's
`momentum cells', whereas for the Einstein metric, the conformal infinity
(which in K\"ahler terms is the zero-set of $s_h$) separates the quotients of
$\HH2$ from the quotients of $\HH{1,1}$. Also some of the self-dual K\"ahler
quotients of $S^5$ will have associated Einstein metrics which are
scalar-flat, while some of the SDE quotients of $\HP2$, $\HH2$ and $\HH{1,1}$
will be conformally flat.

One way to relate the classifications is to observe that every element of
$\su(1,3)$ has a spacelike eigenvector, and some of them (the `elliptic'
elements) have a timelike eigenvector too.  Since $\PSU(1,3)$ acts
transitively on the spacelike or timelike lines, we can fix one of each and
conjugate any element of $\su(1,3)$ into $\un(1,2)$, and the elliptic elements
into $\un(3)$.  On the other hand all adjoint orbits in $\symp(3)$ are
represented by elements of $\un(3)$, and the same is true for $\symp(1,2)$,
since we have given representatives in $\un(1,2)$ in Definition~\ref{def2}.

\begin{remark} There is a rather beautiful Hermitian/quaternionic real form
of the classical Klein correspondence that allows us to make the
identification of adjoint orbits more natural. Recall that there is a special
isomorphism between $\so(6,\C)$ and $\sgl(4,\C)$: $\C^4$ is the spin
representation of $\so(6,\C)$, or, more straightforwardly, $\sgl(4,\C)$ acts
on $\Lambda^2\C^4$ (via $A\cdot u\wedge v = A(u)\wedge v+u\wedge A(v)$)
preserving a complex bilinear form $g_c$ given by the contraction of
$(\alpha,\beta)\mapsto\alpha\wedge\beta$ with the volume element. This
isomorphism underlies the Klein correspondence:
\begin{itemize}
\item lines in $P(\C^4)$ correspond bijectively to points on the quadric in
$P(\Lambda^2\C^4)$ ($P(U)$ corresponds to null line $\Lambda^2U$);
\item points in $P(\C^4)$ correspond bijectively to $\alpha$-planes in the
quadric ($[u]$ corresponds to projectivization of the maximal totally null
subspace $\{u\wedge v:v\in \C^4\}$);
\item planes in $P(\C^4)$ correspond bijectively to $\beta$-planes in the
quadric ($P(W)$ corresponds to the projectivization of the maximal totally
null subspace $\Lambda^2W$).
\end{itemize}
Now $\su(1,3)$ is the real form of $\sgl(4,\C)$ preserving a Hermitian metric
$(.,.)$ of signature $(1,3)$. Consider now the Hodge star operator on
$\Lambda^2\C^{1,3}$ defined by $({*\alpha})\wedge\beta=(\alpha,\beta)\vol$.
For this to make sense, we must take the Hermitian metric to be anti-linear in
$\alpha$ and thus $*$ anti-commutes with $i$. The signature of the metric
implies that $*^2=-1$, so $j:=*$ defines a quaternionic structure on
$\Lambda^2\C^{1,3}$.  It is convenient to make $\Lambda^2\C^{1,3}$ into a
\emph{right} quaternionic vector space in this way (thus $k=ij=j\circ i$).

We denote by $\so^*(3,\HQ)$ the subalgebra of $\so(6,\C)$ commuting with $j$:
it is the real form isomorphic to $\su(1,3)$. We can describe it in
quaternionic terms as the Lie algebra of the group of $\HQ$-linear
transformations of $\HQ^3$ preserving an $(i,j,k)$-invariant skew form
$\omega$, and hence also the triple of signature $(6,6)$ symmetric forms
$g_i,g_j,g_k$ defined by $g_i(a,b)=\omega(ai,b)$ and so on. Note that $g_i$ is
$i$-invariant, but is anti-invariant with respect to $j$ and $k$, and
similarly for $g_j$ and $g_k$. Hence the quaternionic definition is related to
the complex one by taking $g_j$ to be the real part of $g_c$ (since $g_c$ is
$i$-bilinear and $j$-invariant).

A spacelike or timelike line in $\C^{1,3}$ defines a maximal totally null
$(\alpha)$ subspace of $\Lambda^2\C^{1,3}$ and its perpendicular hyperplane
defines a complementary maximal totally null $(\beta)$ subspace.  Such a
decomposition is equivalently given by a $g_j$-orthogonal complex structure
$I$ on $\Lambda^2\C^{1,3}$ commuting with the quaternionic structure: the null
subspaces are the $\pm i$ eigenspaces. Note that $g(a,b)=\omega(Ia,b)$ is
therefore an $(i,j,k)$-invariant inner product and it is easy to check that it
is indefinite or definite according to whether the line is spacelike or
timelike.

An element of $\su(1,3)$ belongs to $\un(1,2)$ or $\un(3)$ (i.e., preserves
the spacelike or timelike line) if and only if its action on
$\Lambda^2\C^{1,3}$ commutes with $I$ if and only if it is skew with respect
to $g$. In fact this realizes $\un(1,2)$ and $\un(3)$ as
$\symp(1,2)\cap\so^*(3,\HQ)$ and $\symp(3)\cap\so^*(3,\HQ)$ respectively.

For example, consider the diagonal element
\begin{equation*}
\begin{pmatrix}
ir_0 & 0 & 0 & 0 \\
0  & ir_1 & 0 & 0 \\
0  & 0  & ir_2 & 0 \\
0  & 0  & 0   & ir_3
\end{pmatrix}
\end{equation*}
with $r_0+r_1+r_2+r_3=0$, defined using the standard basis $e_0,e_1,e_2,e_3$
for $\C^{1,3}$ with $e_0$ timelike. Its action on $\Lambda^2\C^{1,3}$ with
respect to the quaternionic basis $e_0\wedge e_1$, $e_0\wedge e_2$, $e_0\wedge
e_3$ is easily computed to be
\begin{equation*}
\begin{pmatrix}
i(r_0+r_1) & 0 & 0 \\
0  & i(r_0+r_2) & 0  \\
0  & 0  & i(r_0+r_3) 
\end{pmatrix},
\end{equation*}
where $i$ acts by left multiplication (we have chosen our quaternionic basis
so that the complex structure $I$ determined by $e_0$ is left multiplication
by $i$).
\end{remark}

Adjoint orbits in $\su(1,3)$ are essentially determined by their
characteristic and minimal polynomials, and Bryant~\cite{MR2002i:53096} gives
his classification in these terms---more precisely, in terms of the
polynomials of the associated Hermitian matrices. If $P_c$ is the
characteristic polynomial and $P_m$ is the minimal polynomial, then the degree
$d$ of $P_c/P_m$ determines the local cohomogeneity of the self-dual K\"ahler
metric as $2-d$.  In the generic, local cohomogeneity two, case Bryant
discusses the classification in detail, which he divides into Cases 1--4.

For reference, we shall give the correspondence between adjoint orbits in
$\su(1,3)$ and $\symp(1,2)$ which relate Bryant's classification to ours.  We
do this by giving the characteristic and minimal polynomials $P_c$
corresponding to the representatives in Definition~\ref{def2}.  These
correspondences are obtained by choosing an element of $\su(1,3)$ with given
$P_c$, $P_m$ and spacelike eigenvector $e_1$, and computing its action on
$\Lambda^2\C^{1,3}$ with quaternionic basis $e_1\wedge e_0$, $e_1\wedge e_2$,
$e_1\wedge e_3$.

We recall that there are exceptional adjoint orbits which do not give SDE and
self-dual K\"ahler metrics which are conformal. For these exceptional orbits,
the self-dual K\"ahler metric is conformal to a scalar-flat SDE metric, while
the quotient SDE metric is the real hyperbolic metric (conformally flat).

We begin with the (local) cohomogeneity two K\"ahler metrics, where
$P_m(t)=P_c(t)$.

\smallbreak (i) $P_c(t)=(t-r_0)(t-r_1)(t-r_2)(t-r_3)$, where $r_0,r_1,r_2,r_3$
are distinct with $r_0+r_1+r_2+r_3=0$. This is Bryant's Case 4, and
corresponds to
\begin{equation*}
T_0(ip_0,ip_1,ip_2)=\ \  
\begin{pmatrix}
ip_0&0&0\\
0&ip_1&0\\
0&0&ip_2
\end{pmatrix}.
\end{equation*}
with $p_0=r_0+r_1$, $p_1=-(r_0+r_2)$, $p_2=-(r_0+r_3)$ (and so $p_i\neq \pm
p_j$ for $i,j$ distinct). The exceptional orbits arise when one of the weights
vanish.

\smallbreak (ii) $P_c(t)=(t-r_1)(t-r_2)(t-r-i\lambda)(t-r+i\lambda)$, where
$r_1,r_2$ are distinct with $r_1+r_2+2r=0$. This is Bryant's Case 1
and corresponds to
\begin{equation*}
T_0(\lambda,ip,iq)=\ \
\begin{pmatrix}
ip &\lambda &0\\
\lambda &ip&0\\
0&0&iq
\end{pmatrix}
\end{equation*}
with $p=r+r_1$, $q=-2r=r_1+r_2$ (and so $p\neq 0$). The exceptional orbits
arise when $q$ vanishes.

\smallbreak (iii) $P_c(t)= (t-r_1)(t-r_2)(t-r)^2$, where $r_1,r_2$ and $r$ are
distinct with $r_1+r_2+2r=0$. This is Bryant's Case 3 and corresponds to
\begin{equation*}
T_1(1,ip,iq)=\ \
\begin{pmatrix}
ip &0 &0\\
0 &ip&0\\
0&0&iq
\end{pmatrix}+
\begin{pmatrix}
i&i&0\\
-i&-i&0\\
0&0& 0
\end{pmatrix}
\end{equation*}
with $p=r+r_1$, $q=-2r=r_1+r_2$ (and so $p\neq 0$ and $p\neq \pm q$). The
exceptional orbits arise when $q=0$.

\smallbreak (iv) $P_c(t)=(t-r_1)(t-r)^3$, where $r_1$ and $r$ are distinct
with $r_1+3r=0$. This is Bryant's Case 2 and corresponds to
\begin{equation*}
T_2(1,ip)=\ \ 
ip\,\bbi_3+
\begin{pmatrix}
0&0&-i\\
0&0&i\\
i&i&0
\end{pmatrix}
\end{equation*}
with $p=r+r_1$ (so that $p\neq 0$).

\smallbreak
We finally consider the cohomogeneity one and homogeneous K\"ahler metrics.

\smallbreak (i) $P_c(t)=(t-r_0)^2(t-r_1)(t-r_2)$ and
$P_m(t)=(t-r_0)(t-r_1)(t-r_2)$, where $r_0,r_1,r_2$ are distinct with
$2r_0+r_1+r_2=0$. These metrics have cohomogeneity one under $\Un(2)$ or
$\Un(1,1)$ according to the signature of the Hermitian metric on the repeated
eigenspace, and correspond to $T_0(ip,\pm ip,iq)$ or $T_0(iq,ip,\pm iq)$ with
$p\neq\pm q$. When $q=0$ we have an exceptional orbit.

Further degenerations give homogeneous metrics:
\begin{itemize}
\item $P_c(t)=(t-r)^3(t+3r)$ and $P_m(t)=(t-r)(t+3r)$ with $r\neq 0$
corresponds to $T_0(ip,ip,ip)$ and the Bergman metric;
\item $P_c(t)=(t-r)^2(t+r)^2$ and $P_m(t)=(t-r)(t+r)$ with $r\neq 0$ is
exceptional: the K\"ahler metric is the product metric on $S^2\times\cH^2$ and
the SDE quotient (by $T_0(0,0,ip)$ or $T_0(ip,0,0)$) is $\cH^4$.
\end{itemize}

\smallbreak (ii) $P_c(t)=(t+r)^2(t-r-i\lambda)(t-r+i\lambda)$ and $P_m(t)=
(t+r)(t-r-i\lambda)(t-r+i\lambda)$. These have cohomogeneity one under
$\Un(2)$ and correspond to the Pedersen metrics $T_0(\lambda,0,ir)$,
apart from exceptional orbits when $r=0$.

\smallbreak (iii) $P_c(t)= (t+r)^2(t-r)^2$ or $(t+r)(t-r)^3$ and
$P_m(t)=(t+r)(t-r)^2$ with $r\neq 0$, correspond to the height one quotients
by $T_1(1,0,iq)$ and $T_1(1,ip,\pm ip)$, which are cohomogeneity one metrics.

A further degeneration gives an exceptional orbit: $P_c(t)=t^4$ and
$P_m(t)=t^2$. The K\"ahler metric in this case is flat, while the SDE quotient
by $T_1(1,0,0)$ is $\cH^4$.

\smallbreak (iv) $P_c(t)=t^4$ and $P_m(t)=t^3$ corresponds to the height two
quotient $T_2(1,0)$. This is an exceptional orbit: the self-dual K\"ahler
metric has cohomogeneity one, but the SDE quotient is $\cH^4$.

\bibliographystyle{amsalpha}

\providecommand{\bysame}{\leavevmode\hbox to3em{\hrulefill}\thinspace}
\providecommand{\MR}{\relax\ifhmode\unskip\space\fi MR }
\providecommand{\MRhref}[2]{%
  \href{http://www.ams.org/mathscinet-getitem?mr=#1}{#2}
}
\providecommand{\href}[2]{#2}

\end{document}